\documentclass[11pt,a4paper]{article}
\usepackage[top=2.5cm, bottom=2.5cm, left=2.5cm, right=2.5cm]{geometry}
\usepackage[
  pagebackref=false,
  colorlinks,
]{hyperref}

\usepackage{amsmath}
\usepackage{amsthm}
\usepackage{amsfonts}
\usepackage{amssymb}
\usepackage[all]{xy}
\usepackage{graphicx}
\usepackage{tikz-cd}
\usepackage{mathrsfs}
\usepackage[utf8]{inputenc}

\usepackage{mathtools}

\usepackage[affil-it]{authblk}

\hypersetup{citecolor=blue}

\theoremstyle{plain}
\newtheorem{theorem}{Theorem}[section]
\newtheorem{lemma}[theorem]{Lemma}
\newtheorem{proposition}[theorem]{Proposition}
\newtheorem{proposition-definition}[theorem]{Proposition-Definition}

\theoremstyle{remark}
\newtheorem{remark}[theorem]{Remark}

\newtheoremstyle{custombold}
  {\topsep}   
  {\topsep}   
  {\itshape}  
  {}          
  {\bfseries} 
  {.}         
  {.5em}      
  {}          

\theoremstyle{custombold}

\newtheorem{innercustomtheorem}{Theorem}
\newenvironment{customtheorem}[1]
  {\renewcommand\theinnercustomtheorem{#1}\innercustomtheorem}
  {\endinnercustomtheorem}

\newtheorem{innercustomcorollary}{Corollary}
\newenvironment{customcorollary}[1]
  {\renewcommand\theinnercustomcorollary{#1}\innercustomcorollary}
  {\endinnercustomcorollary}

\newcommand{\irm}{\mathrm{i}}
\DeclareMathOperator{\Hom}{Hom}

 \DeclareMathOperator{\im}{\mathrm{Im}}
\DeclareMathOperator{\vol}{Vol}
\newcommand{\re}{{\operatorname{Re}\,}}

\usepackage{titlesec}

\title{On the generating series of the degree sequence\date{\vspace{-2\baselineskip}}
}
\author{{\normalfont  Quang-Khai Nguyen}}
\begin{document}
\maketitle
\begin{abstract}
	We study the generating series associated with the degree sequence of a monomial self-map of a projective toric variety. We establish conditions under which this series has its circle of convergence as a natural boundary, and hence is a transcendental, non-holonomic function. In the case of toric surfaces, our results are sharp; moreover, we answer a question of Bell by proving that its reduction modulo 
$p$ is transcendental for all but finitely many prime numbers $p$.
\end{abstract}
\newcommand\blfootnote[1]{%
  \begingroup
  \renewcommand\thefootnote{}\footnote{#1}%
  \addtocounter{footnote}{-1}%
  \endgroup
}
\blfootnote{2020 \emph{Mathematics Subject Classification}: 37F10, 14M25, 11J81.}
\blfootnote{\emph{Key words and phrases}: Dynamical degrees, Toric geometry, Transcendence.}
\section{Introduction}\label{Introduction}
Degree sequences and the associated dynamical degrees play a central role in algebraic dynamics.
They are defined for arbitrary dominant rational maps on normal projective varieties over any field (see~\cite{Dinh-Sibony-2005,Truong-2020,Dang-2020}).
As fundamental birational invariants, dynamical degrees govern key dynamical features of the system and measure its complexity.
Although their general behavior remains largely mysterious (cf.~\cite{BDJ-2020}), they can be explicitly computed and are known to be algebraic in certain situations, for instance, for dominant monomial maps. In this paper, we restrict ourselves to this latter setting.
Dynamical degrees are defined via the asymptotic behavior of degree sequences, which therefore encode finer information and are, in principle, more difficult to study.
Our main objective is to show that, although dynamical degrees are comparatively simple and well understood in this case, the associated degree sequences may already display non-trivial complexity.
We capture this complexity through the study of the corresponding generating series.

 Throughout this paper, we work in the following setting.
Let $X(\Sigma)$ be a projective toric variety of dimension $d\geq2$ over an algebraically closed field of arbitrary characteristic. The variety $X(\Sigma)$ is determined by a lattice $N\cong\mathbb Z^d$, its dual lattice $M$, and a complete fan $\Sigma$ in $N_{\mathbb R}$. Let $\varphi\colon X(\Sigma)\dashrightarrow X(\Sigma)$ be a dominant monomial map. This map corresponds to a homomorphism $M\rightarrow M$ given by a matrix $A\in {M}_d(\mathbb Z)$ of non-zero determinant. Let $0\leq k\leq d$ be an integer, and let $D$ be an ample Cartier divisor on $X(\Sigma)$; then the \emph{$k$-degree sequence} is defined by the intersection numbers $\deg_{D,k}(\varphi^n)=((\varphi^n)^*D^k\cdot D^{d-k})$ for all $n\geq 0$ (we use the convention $\varphi^0=\mathrm{id}$). This sequence does not depend on the linear equivalence class of $D$. The \emph{$k$-dynamical degree} of $\varphi$ is defined as $\lambda_k(\varphi)\coloneqq \lim \deg_{D,k}(\varphi^n)^{{1}/{n}}$, which is equal to the product of the $k$ eigenvalues of largest modulus of $A$, see~\cite{Hasselblatt-Propp-2007,Favre-Wulcan-2012,Lin-2012-BSMF}. We denote $\lambda_k(\varphi)$ by $\lambda_k$, unless otherwise specified. Our main object of interest is the \emph{generating power series}
$$\Delta_{k,\varphi,D}(z)\coloneqq \sum_{n\geq0}\deg_{D,k}(\varphi^n)z^n\in\mathbb Z[[z]],$$ which converges inside the open disk $\mathbb D(0,\lambda_k^{-1})$. 

These series have been studied extensively in~\cite{Abarenkova-Auriac-Boukraa-Maillard-1999,Boukraa-Hassani-Maillard-2003} when $k=1$ and have been proved to be rational in many cases, most notably for regular maps and for birational surface maps. When the base field is $\mathbb C$, the dynamical degrees and the degree sequences are well-understood when the dynamical system $(X(\Sigma),\varphi)$ can be lifted to a simplicial projective toric variety $X(\Tilde{\Sigma})$ and a \emph{$k$-stable} self-map $\Tilde{\varphi}$ in the sense that $\Tilde{\varphi}^{n*}=\Tilde{\varphi}^{*n}$ for the induced pullbacks on $H^{2k}(X(\Tilde{\Sigma}),\mathbb R)$ and for all $n\geq0$, a property first introduced in~\cite{Fornaess-Sibony-1995}. A consequence of  $k$-stability is that the degree sequences satisfy a \emph{linear recurrence relation} (see~\cite[Corollary~2.2]{Diller-Favre-2001} and~\cite[Proposition~6.3]{Lin-Wulcan-2014}), which implies that $\Delta_{k,\varphi,D}$ is a \emph{rational function}. On the other hand, under certain conditions on $A$, it is shown in~\cite{Bedford-Kim-2004,Hasselblatt-Propp-2007,Lin-Wulcan-2014} that the degree sequences do not satisfy any linear recurrence, so that  $\Delta_{k,\varphi,D}$ is not rational.

A natural way to measure the complexity of power series is to place them within the following hierarchy:
\begin{equation}\label{eq:classification}
 \mathrm{Rat} \subset \mathrm{Alg} \subset \mathrm{Hol} \subset \mathrm{DiA} \subset \mathbb{C}[[z]].
\end{equation}
Here, $\mathrm{Rat} = \mathbb{C}(z) \cap \mathbb{C}[[z]]$ denotes the ring of rational power series, $\mathrm{Alg}$ the ring of power series that are \emph{algebraic} over $\mathbb{C}(z)$ (i.e., those satisfying a polynomial equation with coefficients in $\mathbb{C}(z)$), $\mathrm{Hol}$ the ring of \emph{holonomic} power series (that is, those satisfying a linear differential equation with coefficients in $\mathbb{C}(z)$), and $\mathrm{DiA}$ the ring of \emph{differentially algebraic} power series, those satisfying an algebraic differential equation with coefficients in $\mathbb C(z)$. 
Power series that do not belong to $\mathrm{Alg}$ (resp. $\mathrm{DiA}$) are called \emph{transcendental} (resp. \emph{hypertranscendental} or \emph{differentially transcendental}).  
This classification also provides information about the corresponding coefficient sequences. Indeed, a power series is rational if and only if its coefficients satisfy a linear recurrence with constant coefficients, while it is holonomic if and only if its coefficients satisfy a linear recurrence relation with polynomial coefficients (see~\cite[Chapter~4]{Stanley-Enumerative-Combinatorics-2011} and~\cite[Chapter~6]{Stanley-Fomin-vol2-1999}). Thus, more involved functional equations give rise to more complex coefficient sequences. 
When restricted to power series with algebraic coefficients, the rings $\mathrm{Hol}$ and $\mathrm{DiA}$ are countable, and they play for arithmetic holomorphic functions a role analogous to that of the ring of periods among complex numbers~\cite{Kontsevich-Zagier-2001}. The classification~\eqref{eq:classification} is also central in enumerative combinatorics, where the nature of a generating series reflects structural properties of the objects being counted (cf.~\cite{Bousquet-Melou-2006}). 
It has long been a classical problem to determine to which ring a given series arising from number theory, combinatorics, complex dynamics, or physics belongs. 
Hilbert~\cite{Hilbert-1902}   already observed that many important number-theoretic functions, such as the Riemann zeta function, do not belong to $\mathrm{DiA}$. For further results in this direction, we refer to~\cite{Stanley-1980,Becker-Topfer-94,Stanley-Fomin-vol2-1999,Adamczewski-Dreyfus-Hardouin-2021} and the references therein.

 A distinctive feature of holonomic convergent power series is that they admit analytic continuation to all but finitely many points of 
$\mathbb C$ (see, for instance,~\cite[Theorems~2.1 and 4(a)]{Stanley-1980}). A natural strategy to prove that a power series is not algebraic or holonomic is therefore to show that it has infinitely many singularities. In the case of the series $\Delta_{k,\varphi,D}(z)$, we investigate its analytic properties in greater detail and establish a stronger result. Namely, we will prove that it admits the circle $\mathcal{C}(0,{\lambda_k}^{-1})$ as a \emph{natural boundary}, i.e., it cannot be analytically continued at any point on $\mathcal{C}(0,{\lambda_k}^{-1})$. Our main theorem, Theorem~\ref{natural boundary of k-degree sequence}, may be viewed as a generalization of~\cite[Th\'eor\`eme principal]{Favre-2003},~\cite[Theorem~1.1]{Bedford-Kim-2004},~\cite[Proposition~7.6]{Hasselblatt-Propp-2007}, and~\cite[Theorem~C]{Lin-Wulcan-2014}. The existence of a natural boundary shows that even in relatively simple cases where the dynamical degrees can be computed explicitly, the associated degree sequences may nevertheless exhibit some complexity.

\begin{customtheorem}{A}
\label{natural boundary of k-degree sequence}
     Let $\varphi\colon X(\Sigma)\dashrightarrow X(\Sigma)$ be a dominant monomial self-map of a projective toric variety of dimension $d\geq2$, which corresponds to a matrix $A \in M_d(\mathbb Z)$ with $\det A\neq0$.  Let $D$ be an ample Cartier divisor on $X(\Sigma)$.  Assume that the eigenvalues of $A$ are ordered as $|\mu_1|\geq|\mu_2|\geq|\mu_3|\geq\ldots\geq| \mu_d|$.
     Let $1\leq k\leq d-1$, and assume further that 
     \begin{itemize}
         \item $|\mu_{k-1}|>|\mu_k|=|\mu_{k+1}|>|\mu_{k+2}|$,
         \item $\mu_{k}$ and $\mu_{k+1}$ are complex conjugates, and $\mu_{k+1}/\mu_{k}$ is not a root of unity.
     \end{itemize}
     Then the series $\Delta_{k,\varphi,D}(z)$ has the circle $\mathcal{C}(0,{\lambda_k}^{-1})$ as a natural boundary. In particular, the series $\Delta_{k,\varphi,D}(z)$ is transcendental and non-holonomic over $\mathbb{C}(z)$. 
\end{customtheorem}

\begin{remark}\label{remark1}
There are several settings in which it has been shown that a power series is either rational or has its circle of convergence as a natural boundary. See, for instance,~\cite{Bell-Coons-Rowland-2013} in the context of Mahler equations, and~\cite{Nguyen-Bell-Gunn-Saunders-2023,Boudreau-Holmes-Nguyen-2025} in the context of Artin--Mazur zeta functions arising from dynamics. It would be interesting to know whether a similar dichotomy also holds in our setting; in this direction, see Theorem~\ref{dichotomy of the generating power series for toric surfaces}.
We stress, however, that proving the existence of a natural boundary does not imply that the corresponding function is differentially transcendental. Establishing the latter property usually requires more sophisticated algebraic tools and is often related to the existence of a functional equation (not of differential type) satisfied by the function in question (see~\cite{Adamczewski-Dreyfus-Hardouin-2021}). It would therefore be interesting to determine whether the generating series in Theorem~\ref{natural boundary of k-degree sequence} are differentially transcendental.
\end{remark}

\begin{remark}\label{Artin--Mazur style zeta function}
One can consider other variants of the generating series, such as the Artin--Mazur style zeta function $$\zeta_{k,\varphi,D}(z)\coloneqq \mathrm{exp}\Bigg(\sum_{n\geq1}\frac{\deg_{D,k}(\varphi^n)}{n}z^n\Bigg)$$ which also converges in $\mathbb{D}(0,\lambda_k^{-1})$. Observe that
$z\frac{d}{dz}\log \zeta_{k,\varphi,D}(z)=\Delta_{k,\varphi,D}(z)-(D^d)$. Therefore $\zeta_{k,\varphi,D}(z)$ also has $\mathcal C(0,\lambda_k^{-1})$ as a natural boundary.

The Ces\`aro mean of the sequence $\deg_{D,k}(\varphi^n)$ is given by 
	$$\sum_{n\geq1}\frac{\sum_{0\leq m\leq n-1}\deg_{D,k}(\varphi^m)}{n}z^n.$$
	A computation similar to Theorem~\ref{natural boundary of k-degree sequence} implies that under the same assumption, it also admits $\mathcal C(0,\lambda_k^{-1})$ as a natural boundary.
\end{remark}

The proof of Theorem~\ref{natural boundary of k-degree sequence} proceeds as follows. First, we use the standard dictionary relating intersection theory in toric varieties and mixed volumes to compute the degrees, following~\cite{Huber-Sturmfels-1995}. It expresses the degree sequence as the values of
a piecewise linear function at the sequence of powers of an auxiliary matrix. Under our assumption on the eigenvalues of $A$, we are then able to reduce the whole problem to a two-dimensional situation where the auxiliary matrix
becomes conformal. At this stage, we need to control a series of the form $\sum h(\chi^n)z^n$ where $h$ is piecewise linear on the complex plane, and $\chi$ is a complex number whose argument is incommensurable with $2\pi$. The key argument is to consider the Fourier series of the restriction of $h$ to the unit circle, and to exploit the discontinuity of the derivative of $h$ to infer that infinitely many of its Fourier coefficients
are non-zero. To conclude that the generating series has a natural boundary, we use in a crucial way that
these coefficients are non-zero along some arithmetic progression, which follows
from the Skolem--Mahler--Lech theorem.

In~\cite{BDJ-2020}, when the characteristic of the base field is different from 2, the authors construct the first example of a dominant rational self-map $f_\mu$ of $\mathbb P^2$ with transcendental dynamical degree $\lambda_1(f_\mu)$. The map $f_\mu$ is of the form $f_\mu=g\circ \varphi_\mu$, where   $$g\colon \mathbb P^2\dashrightarrow \mathbb P^2, [x_0 : x_1 : x_2] \mapsto [x_0(x_1 + x_2 - x_0) : x_1(x_2 + x_0 - x_1) : x_2(x_0 + x_1 - x_2)],$$
is a birational involution of $\mathbb P^2$, $\varphi_\mu$ is a monomial map associated with \[A_\mu=\begin{pmatrix}
\re \mu & -\im\mu \\ 
 \im \mu & \re \mu
\end{pmatrix}\] where $\mu\in\mathbb Z[\irm]$ is a Gaussian integer whose argument is incommensurable
with $2\pi$. Recall from \emph{loc.~cit.} that $\lambda_{1}(f_\mu)>\lambda_1(\varphi_\mu)=\lambda_2(\varphi_\mu)^{1/2}=\lambda_2(f_\mu)^{1/2}$ and we have the identity:
$$ (1+\Delta_{1,f_\mu,\mathcal O(1)}(z))(2-\Delta_{1,\varphi_\mu,\mathcal O(1)}(z))=2$$
for  $|z|<\lambda_{1}(f_\mu)^{-1}$, see~Formula (2.2) and Proposition~2.8 in \emph{loc.~cit.}. It is straightforward to deduce the following corollary by applying Theorem~\ref{natural boundary of k-degree sequence} and~\cite[Main Theorem]{Boucksom-Favre-Jonsson-2008}.  

\begin{customcorollary}{B}\label{the example of BDJ} For any Gaussian integer $\mu$ whose argument is incommensurable with $2\pi$, the generating series $\Delta_{1,f_\mu,\mathcal{O}(1)}(z)$ of the map $f_\mu = g \circ \varphi_\mu$ is meromorphic in the disk $\mathbb{D}\bigl(0, \lambda_2(f_\mu)^{-1/2}\bigr)$, where it has a simple pole at $\lambda_1(f_\mu)^{-1}$ and no other poles. In addition, it admits the circle $\mathcal{C}(0, \lambda_{2}(f_\mu)^{-1/2})$ as a natural boundary.
\end{customcorollary}

Let $\mathfrak M$ denote the category of all normal projective varieties $X$ over arbitrary fields $K$, together with all ample Cartier divisors $D$, and all dominant rational self-maps $f$. It is shown in~\cite{Urech-2018} that the set of degree sequences attached to elements of $\mathfrak M$ is countable\footnote{Although the result in \emph{loc.~cit.} is stated only for smooth varieties, its proof extends verbatim to the normal case.}. In view of the examples of transcendental dynamical degrees~\cite{BDJ-2020,BDJK-2024,Sugimoto-2025}, it is thus natural to ask whether $$\mathrm{tr.deg}_{\overline{\mathbb Q}}(\lambda_k(f):(X,K,D,f)\in \mathfrak M)=\infty \;\;\;\; \forall k\geq1.$$ However, even the weaker question $\mathrm{tr.deg}_{\overline{\mathbb Q}}(\lambda_k(f):(X,K,D,f)\in \mathfrak M)\geq2$ remains open. On the other hand, at the level of functions, Theorem~\ref{natural boundary of k-degree sequence} allows us to deduce the following.

\begin{customcorollary}{C}\label{algebraically independence of generating functions}
For every $k\geq1$, we have $\mathrm{tr.deg}_{\mathbb C(z)}(\Delta_{k,f,D}(z):(X,K,D,f)\in \mathfrak M)=\infty.$
\end{customcorollary}

The proof of Corollary~\ref{algebraically independence of generating functions} in fact establishes a stronger statement: for any monomial map $\varphi$ satisfying the assumptions of Theorem~\ref{natural boundary of k-degree sequence}, the series $\Delta_{k,\varphi^\ell,D}(z), \ell\geq1,$ are algebraically independent over $\mathbb C(z)$.

There is a parallel story for the degree sequences and dynamical degrees in the context of \emph{monomial correspondences}~\cite{Dang-Ramadas-2021}. Combining our results with those in \emph{loc.~cit.}, it is straightforward to establish a version of Theorem~\ref{natural boundary of k-degree sequence} and its corollaries in this setting. To keep the paper reasonably short, we refrain from stating it.

In the two-dimensional case, Theorem~\ref{natural boundary of k-degree sequence} covers all cases, and we obtain:
\begin{customtheorem}{D}\label{dichotomy of the generating power series for toric surfaces}    Let $\varphi\colon X(\Sigma)\dashrightarrow X(\Sigma)$ be a dominant monomial self-map of a projective toric surface. Let $A\in M_2(\mathbb Z)$ be the matrix associated with $\varphi$, and let $\mu_1$, $\mu_2$ be its eigenvalues. Let $D$ be an ample Cartier divisor on $X(\Sigma)$.
\begin{enumerate}
    \item If $\mu_1$ and $\mu_2$ are complex conjugates and $\mu_1/\mu_2$ is not a root of unity, then  $\Delta_{1,\varphi,D}$ admits $\mathcal{C}(0,\lambda_1^{-1})$ as a natural boundary. Furthermore, for any prime number $p$ sufficiently large, the reduction modulo $p$  
    $$\Delta_{1,\varphi,D}(z)\bmod p\coloneqq \sum_{n\geq0}(\deg_{D,1}(\varphi^n)\bmod p)z^n\in\mathbb F_p[[z]]$$
    is transcendental over $\mathbb F_p(z)$.

    \item Otherwise, $\Delta_{1,\varphi,D}$ is a rational function.
\end{enumerate}
\end{customtheorem}

The proof of the dichotomy in Theorem~\ref{dichotomy of the generating power series for toric surfaces} relies on~\cite[Th\'eor\`eme~principal]{Favre-2003} 
(see also~\cite[Theorem~C']{Jonsson-Wulcan-2011}), which provides a necessary and sufficient condition for the existence of algebraically stable models for monomial self-maps of toric surfaces.
Moreover, the transcendence of the reductions modulo 
$p$ of the generating series answers a question posed by Bell (see~\cite{Bell-talk-NTWS}). Our argument in this direction relies on Christol's theorem from automata theory.

\begin{remark}
We stress that the study of reductions modulo $p$ of power series in $\mathbb Z[[x]]$ is an active topic, another way to measure the complexity of the coefficient sequence. When such a power series is algebraic, all its reductions modulo $p$ are necessarily algebraic. The converse, however, is false: there exist transcendental power series whose reductions are algebraic for every prime $p$. In particular, this intriguing phenomenon occurs for important families of holonomic power series arising from geometry (such as diagonals of algebraic power series~\cite{Deligne-1984} and some hypergeometric series~\cite{Vargas-Montoya-2021}).  It has proved useful in questions of functional transcendence and has also led to challenging open problems (see, for instance,~\cite{Adamczewski-Bell-2013,Adamczewski-Bell-Delaygue-2019,Adamczewski-Bostan-Caruso-2023,Caruso-Furnsinn-Vargas-Montoya-2025}).
\end{remark}

\begin{remark}
	Note that by Remark~\ref{Artin--Mazur style zeta function}, the same dichotomy holds for the Artin--Mazur zeta function $\zeta_{1,\varphi,D}(z)$ in the toric surface case. Similar results have been proved using the P\'olya--Carlson theorem for other kinds of Artin--Mazur zeta functions, see~\cite[Corollary~10.4]{Milnor-Thurston-1988}
 and~\cite[Theorem~4.2]{Nguyen-Bell-Gunn-Saunders-2023} for examples. However, the latter papers put conditions on the coefficients of the series that are not satisfied in our case.
\end{remark}

Finally, we prove an intriguing rigidity result for the generating series of monomial surface maps, which follows from Theorem~\ref{dichotomy of the generating power series for toric surfaces}. We refer to~\cite{Ji-Xie-2023} and the references therein for deep rigidity results in holomorphic dynamics in dimension 1.

To state our result, we need an auxiliary notion. A matrix $A \in M_2(\mathbb Z)$ is said to be \emph{absolutely irreducible over $\mathbb Q$} if there exists no $A^n$-invariant line in $\mathbb Q^2$ for all $n\geq1$, or equivalently, $A^n$ has no rational eigenvalues. A monomial self-map $\varphi$ of a toric surface is said to be \emph{absolutely irreducible} over $\mathbb Q$ if its associated matrix is absolutely irreducible over $\mathbb Q$. 

\begin{customtheorem}{E}\label{when two degree sequences coincide}
    Let $\varphi$ and $\varphi'$ be dominant monomial self-maps of a projective toric surface $X(\Sigma)$ that are absolutely irreducible over $\mathbb Q$. Let $D$ and $D'$ be ample Cartier divisors on $X(\Sigma)$. Assume that $\varphi$, with respect to $D$, and $\varphi'$, with respect to $D'$, have the same $1$-degree sequence, i.e., $$\Delta_{1,\varphi,D}=\Delta_{1,\varphi',D'}.$$
    Then there exists $u\in\{1,2,3,4,6\}$ such that $\varphi^u$ and $\varphi'^u$ are semi-conjugate, that is, $f\circ \varphi^u=\varphi'^u\circ f $ for some dominant monomial map $f:X(\Sigma)\dashrightarrow X(\Sigma)$.
\end{customtheorem}

The difficulty in proving Theorem~\ref{when two degree sequences coincide} is when the series are not rational. In this case, we need to use an explicit expression of $\Delta_{1,\varphi,D}$ and $\Delta_{1,\varphi',D'}$ to compare the coefficients of both sides. One can ask whether $\Delta_{1,f,D} = \Delta_{1,f,D'}$ implies $D= D'$, but it is unclear at this stage how to recover the line bundle from the generating series.

The paper is organized as follows. In Section 2, we cover basic preliminaries on toric geometry. In Section 4, we aim to present the proof of Theorem~\ref{natural boundary of k-degree sequence} and Corollary~\ref{algebraically independence of generating functions}. Finally, we restrict our attention to the surface case in Section 5.

\subsection*{Notations}

Throughout this paper, we use $\mathbb N=\{0,1,2,\ldots\}$. We denote $\arg z\in [0,2\pi)$ for the argument of $z\in\mathbb C$, and $\{x\}\in[0,1)$ for the fractional part of $x\in\mathbb R$. For two functions $f,g\colon\mathbb N\rightarrow \mathbb R$, we write $f\asymp g$ when there exist constants $c_1,c_2>0$ such that $c_1f(n)\leq g(n)\leq c_2f(n)$ for all $n$ sufficiently large.

 \subsection*{Acknowledgment}
I am grateful to my advisors Boris Adamczewski and Charles Favre for their suggestions, helpful discussions, and constant encouragement. I would like to thank the anonymous referee for valuable suggestions that helped improve this article.

This project has received funding from the European Union’s MSCA-Horizon Europe, grant agreement No. 101126554. 

\subsection*{Disclaimer}
Co-Funded by the European Union. Views and opinions expressed are, however, those of the author only and do not necessarily reflect those of the European Union. Neither the European Union nor the granting authority can be held responsible for them.

\noindent \includegraphics[height = 1cm]{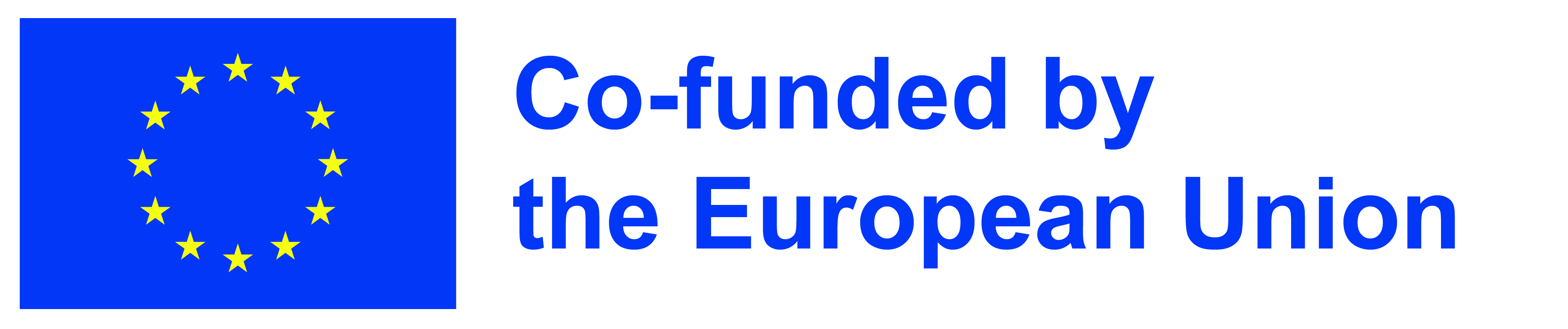}
 \subsection*{Conflict of interest}
There is no conflict of interest.
\subsection*{Data Availability Statement}
This article does not involve any datasets.
\section{Toric geometry}

In this section, we recall some material from toric geometry and aim to interpret the degree sequences using mixed volumes, following~\cite{Oda-toric-geometry}.

A \emph{toric variety} is a normal variety that contains an open dense torus, whose natural multiplicative action on itself extends to the whole variety.
\subsection{Toric varieties}

Let $d\geq 1$. Let $N$ be a free abelian group of rank $d$, and let $M=\Hom_\mathbb Z(N,\mathbb Z)$ be its dual lattice. We denote by $N_\mathbb R\coloneqq N\otimes\mathbb R\cong \mathbb R^d$ and $M_\mathbb R\coloneqq M\otimes\mathbb R$ the associated vector spaces equipped with the natural pairing $\langle\cdot,\cdot\rangle$. We fix a basis of $N$ and the dual basis for $M$.

A \emph{polyhedral convex cone} in $N_\mathbb R$ is a set of the form $\sigma=\Big\{\sum_{v\in S}\lambda_vv:\lambda_v\geq0\Big\}$ for a finite subset $S\subset N_\mathbb R$. A polyhedral convex cone $\sigma$ is said to be \emph{strongly convex} if $\sigma\cap(-\sigma)=\{0\}$, and said to be \emph{rational} if $\sigma$ is generated by a finite subset $S$ of $N$. For the rest of this section, by cones we mean strongly convex rational polyhedral cones. The \emph{dimension} $\dim\sigma$ of $\sigma$ is the dimension of the vector space generated by $\sigma$. One-dimensional cones are also called \emph{rays}. A cone is \emph{simplicial} if it is generated by linearly independent vectors. A cone is \emph{regular} if it is generated by a set which is part of a basis of $N$. 

A \emph{fan} $\Sigma$ in $N_\mathbb R$ is a finite collection of cones such that every face of a cone in $\Sigma$ is also a cone in $\Sigma$, and the intersection of two cones in $\Sigma$ is a face of each. The \emph{support} of $\Sigma$ is defined as 
$|\Sigma|\coloneqq \bigcup_{\sigma\in\Sigma}\sigma$. We denote by $\Sigma(1)$ the set of rays in $\Sigma$. A fan is said to be \emph{simplicial} (resp. \emph{regular})  if all of its cones are simplicial (resp. regular). A fan $\Sigma$ is said to be \emph{complete} if $|\Sigma|=N_\mathbb R$. 
 A fan $\Tilde{\Sigma}$ is a \emph{refinement} of $\Sigma$ if every cone in $\Tilde{\Sigma}$ is included in a cone of $\Sigma$.
 
Each fan $\Sigma$ determines a toric variety $X(\Sigma)$ of dimension $d$ with a dense torus $\mathbb G_m^d$, and all toric varieties arise in this way. The variety $X(\Sigma)$ is complete (resp. regular) if and only if $\Sigma$ is complete (resp. regular). We say that $X(\Sigma)$ is simplicial if $\Sigma$ is simplicial.
   
   We say that $\Sigma$ is \emph{projective} if $X(\Sigma)$ is projective.
\subsection{Cartier divisors and piecewise linear functions} 

Let $\Sigma$ be a complete fan in $N$. On the toric variety $X(\Sigma)$, the action of the dense torus extends to an action on Cartier divisors (resp. Weil divisors). Those invariant under such an action are said to be \emph{$\mathbb G_m^d$-invariant}. Every Cartier divisor (resp. Weil divisor) is linearly equivalent to a $\mathbb G_m^d$-invariant one.

It turns out that $\mathbb G_m^d$-invariant Cartier divisors on $X(\Sigma)$ can be described in terms of support functions as follows. Let $\mathrm{SF}(\Sigma)$ be the set of all \emph{support functions with respect to $\Sigma$}, that is, for every cone $\sigma\in\Sigma$, there exists $m(\sigma)\in M$ such that $h|_\sigma=m(\sigma)$.  For each ray $\rho\in\Sigma(1)$, the associated \emph{ray generator} is the unique primitive lattice point $v_\rho$ generating  $\rho=\mathbb R_{\geq0}v_\rho$. Each ray $\rho$ determines a toric Weil divisor denoted by $D_\rho$. Every $h\in\mathrm{SF}(\Sigma)$ determines a $\mathbb G_m^d$-invariant Cartier divisor $D(h)\coloneqq \sum_{\rho\in \Sigma(1)}h(v_\rho)D_\rho$. Conversely, every $\mathbb G_m^d$-invariant Cartier divisor $D$ can be regarded as a Weil divisor of the form $D=\sum_{\rho\in\Sigma(1)} a_\rho D_\rho$ with $a_\rho\in\mathbb Z$, which defines a function $h_D\in\mathrm{SF}(\Sigma)$ by $h_D(v_\rho)=a_\rho$ for all $\rho\in\Sigma(1)$. This gives a one-to-one correspondence between $\mathrm{SF}(\Sigma)$ and  $\mathbb G_m^d$-invariant Cartier divisors. A $\mathbb G_m^d$-invariant Cartier divisor $D$ is ample if and only if $h_D$ is \emph{strictly convex}, i.e., it is a convex function and defined by different $m(\sigma)$ for different $d$-dimensional cones $\sigma$ in $\Sigma$. Thus $\Sigma$ is \emph{projective} if and only if there is a strictly convex function in $\mathrm{SF}(\Sigma)$.
 
A \emph{lattice polytope} in $M_\mathbb R$ is the convex hull of finitely many points in $M$. A function $h\in \mathrm{SF}(\Sigma)$ determines a  (non-empty) lattice polytope  $P(h)\coloneqq \{m\in M_\mathbb R:m\leq h\}\subset M_\mathbb R$. It follows that a $\mathbb G_m^d$-invariant Cartier divisor $D$  defines a lattice polytope $P_D\coloneqq P(h_D)$. Two  $\mathbb G_m^d$-invariant Cartier divisors $D_1$ and $D_2$ are linearly equivalent if and only if $h_{D_1}-h_{D_2}$ is linear on $N_\mathbb R$,  or equivalently, $P_{D_1}$ is a translate of $P_{D_2}$. We note that if $h\in \mathrm{SF}(\Sigma)$ is strictly convex, then $P(h)$ is a lattice polytope with non-empty
interior, see~\cite[Lemma~1.7.12]{Schneider-book-2013}. 

\subsection{Monomial maps}

 For a homomorphism $A\colon M\rightarrow M$, we denote by $A$ the induced map $M_\mathbb R\rightarrow M_\mathbb R$, and by $\check{A}$ the dual maps $N\rightarrow N$ and $N_\mathbb R\rightarrow N_\mathbb R$.
 
Let $\Sigma$ and $\Sigma'$ be complete fans in $N$. A rational map $\varphi\colon X(\Sigma')\dashrightarrow X(\Sigma)$ is said to be  \emph{monomial} if it is equivariant under the action of the torus $\mathbb G_m^d$. Every group homomorphism $A\colon M\rightarrow M$ with $A=(a_{ij})\in M_d(\mathbb Z)$ gives rise to the group homomorphism $\check{A}\colon N\rightarrow N$, which in turn determines a monomial map
$\varphi\colon\mathbb G_m^d\rightarrow \mathbb G_m^d,(z_1,\ldots,z_d)\mapsto(z_1^{a_{11}}\cdots z_d^{a_{d1}},\ldots,z_1^{a_{1d}}\cdots z_d^{a_{dd}})$
where $z_1,\ldots,z_d$ are coordinates on $\mathbb G_m^d$, hence a monomial map $X(\Sigma')\dashrightarrow X(\Sigma)$. Conversely, any monomial map $\varphi\colon X(\Sigma')\dashrightarrow X(\Sigma)$ arises in this way. The map $\varphi\colon X(\Sigma')\dashrightarrow  X(\Sigma)$ is dominant if and only if  $\det(A)\neq0.$ In this case, its topological degree is equal to $|\det(A)|$.

 When $\varphi\colon X(\Sigma') \rightarrow X(\Sigma)$ is a regular monomial map and $D$ is a $\mathbb G_m^d$-invariant Cartier divisor on $X(\Sigma)$ corresponding to $h_D\in \mathrm{SF}(\Sigma)$ and the polytope $P_D\subset M_\mathbb R$, we have $\varphi^*D=\sum_{\rho'\in\Sigma'(1)}h_D(\check{A}(v_{\rho'}))D_{\rho'}$, which is a $\mathbb G_m^d$-invariant Cartier divisor on $X(\Sigma')$. Further, $\varphi^*D$ corresponds to $h_{\varphi^*D}=h_D\circ \check{A}\in\mathrm{SF}(\Sigma')$ and $P(h_{\varphi^*D})=AP(h_D)\subset M_\mathbb R$.

\subsection{Mixed volumes and intersection numbers}

Let $K_1,\ldots,K_d$ be convex bodies in $M_\mathbb R$. A theorem of Minkowski and Steiner says that the map given by the \emph{Minkowski sum} $$\mathbb R^d_{\geq0}\rightarrow \mathbb R_{\geq0},(r_1,\ldots,r_d)\mapsto\vol(r_1K_1+\cdots+r_dK_d)$$is a homogeneous polynomial of total degree $d$. The \emph{mixed volume} $\vol(K_1,\ldots,K_d)$ of $K_1,\ldots,K_d$ is the non-negative number defined as the coefficient of the monomial $r_1\cdots r_d$, divided by $d!$. 

 By~\cite[p.~79]{Oda-toric-geometry}, if $D_1,\ldots,D_d$ are $\mathbb G_m^d$-invariant Cartier divisors on a complete toric variety $X(\Sigma)$, then the intersection product can be interpreted by mixed volumes as $$(D_1\cdots D_d)=d!\vol(P_{D_1},\ldots,P_{D_d})\in\mathbb Z.$$
 Here, we recall that $P_{D_1},\ldots,P_{D_d}$ are the polytopes associated with $D_1,\ldots,D_d$, respectively. 

\subsection{Degree sequences}

Let $X(\Sigma)$ be a projective toric variety defined by a fan $\Sigma$ in a lattice $N$ of rank $d\geq1$. Let $\varphi\colon X(\Sigma)\dashrightarrow X(\Sigma)$ be a dominant monomial map associated with a homomorphism $A\colon M\rightarrow M$, and let $D$ be an ample Cartier divisor on $X(\Sigma)$.  

For $0\leq k\leq d$, the $k$-\emph{degree} $\deg_{D,k}(\varphi)$ of $\varphi$ with respect to $D$ is defined as follows. We can find a regular refinement $\Tilde{\Sigma}$ of $\Sigma$, which gives rise to a regular modification $\pi\colon X(\Tilde{\Sigma})\rightarrow X(\Sigma)$ induced by the identity map on $N$, such that the composition $\Tilde{\varphi} = \varphi\circ\pi \colon X(\Tilde{\Sigma}) \rightarrow X(\Sigma)$ is a regular monomial map. We define $\deg_{D,k}(\varphi)$ as the intersection number $((\Tilde{\varphi}^*D)^k\cdot(\pi^*D)^{d-k})$. This definition is independent of the choice of $\Tilde{\Sigma}$, and we simply write $\deg_{D,k}(\varphi) = (\varphi^*D^k \cdot D^{d-k})$.

Since every Cartier divisor is linearly equivalent to a $\mathbb G_m^d$-invariant one, we may assume that $D$ is $\mathbb G_m^d$-invariant. The {$k$-degree sequence}  of $\varphi$ can be rewritten as $$\deg_{D,k}(\varphi^n)=d!\vol(A^nP_D[k],P_D[d-k])\in\mathbb Z,$$
see~\cite[Proposition~4.1]{Favre-Wulcan-2012}. Here, the notation $Q[\ell]$ stands for $\ell$ times repetition of $Q$.

\section{Some auxiliary results on piecewise linear functions}

In this section, we present two auxiliary results for real-valued piecewise linear functions that will be useful later. 

A piecewise linear function on a finite-dimensional real vector space $V$ is a continuous 1-homogeneous function $h\colon V \rightarrow \mathbb R$ such that there exists a complete (not necessarily rational) fan $\Sigma$ of $V$ and for each cone $\sigma \in \Sigma$ a linear function $h_\sigma\colon V\rightarrow\mathbb R$ such that
$h = h_\sigma$ on the cone $\sigma$. 

\begin{lemma}\label{piecewise linear functions are Lipschitz}
   Every piecewise linear function is Lipschitz.
\end{lemma}
\begin{proof}
    Let $h\colon V\rightarrow\mathbb R$ be a piecewise linear function with respect to some fan. We fix a norm $|\cdot|$ on $V$. We can find some $C>0$  sufficiently large so that if $v,w\in V$ lie in the same cone, then $|h(v)-h(w)|\leq C|v-w|$. If $v$ and $w$ are not in the same cone, we can join them by a segment $v=v_0,v_1,\ldots,v_n=w$ where $v_{i}$, $v_{i+1}$  stay in the same cone for every $i=0,\ldots,n-1$. Therefore $$|h(v) - h(w)|\leq \sum_{i=0}^{n-1}|h(v_i)-h(v_{i+1})| \leq C \sum_{i=0}^{n-1} |v_i-v_{i+1}| = C |v-w|$$ as desired.
\end{proof}

\begin{proposition}\label{Fourier coefficients of continuous piecewise linear functions}
     Let $g\colon \mathbb R^2\rightarrow \mathbb R$ be a piecewise linear function. Let $\chi$ be some non-zero complex number. We denote by $\theta$ its argument. Then the power series $\sum_{n\geq0}g(\mathrm{Re}(\chi^n),\mathrm{Im}(\chi^n))z^n$ converges inside the disk $\mathbb D(0,|\chi|^{-1})$ and can be written as 
     $$\sum_{m\in\mathbb Z}\frac{a_m}{1-e^
    {\irm m\theta}|\chi|z}$$
     where $a_m\in\mathbb C$ such that $\sum_{m\in\mathbb Z}|a_m|<\infty$. Further, there is a linear recurrence sequence $\{d_m:m\in\mathbb Z\},$ such that $a_m=\frac{d_m}{m^2-1}$ for all $m\neq\pm1$. \end{proposition}
\begin{proof} 
  Since $g$ is piecewise linear, we have $|g(\mathrm{Re}(\chi^n),\mathrm{Im}(\chi^n))|= O(|\chi|^n)$ which yields that $\sum_{n\geq0}g(\mathrm{Re}(\chi^n),\mathrm{Im}(\chi^n))z^n$ converges inside the disk $\mathbb D(0,|\chi|^{-1})$. 
  
  It is convenient to identify $\mathbb R^2$ with $\mathbb C$ to make use of the conformal structure on $\mathbb C$, and view $g$ as a  1-homogeneous continuous function on $\mathbb C$. Then  $$\sum_{n\geq0}g(\mathrm{Re}(\chi^n),\mathrm{Im}(\chi^n))z^n=\sum_{n\geq0}|\chi|^ng(e^{\irm n\theta})z^n$$
   inside the disk $\mathbb D(0,|\chi|^{-1}).$

We decompose the function $g|_{\mathbb S^1}$ in the Hilbertian basis consisting of $\{e^{\irm mx}:m\in\mathbb Z\}$ in $L^2(\mathbb S^1)$. Since $g|_{\mathbb S^1}$ is a piecewise linear function on $\mathbb S^1$, it is Lipschitz by Lemma~\ref{piecewise linear functions are Lipschitz}. Therefore, its Fourier series converges uniformly to $g|_{\mathbb S^1}$. In addition, since $g|_{\mathbb S^1}$ is piecewise linear, $g|_{\mathbb S^1}$ belongs to the Sobolev space $W^{1,2}(\mathbb S^1)\coloneqq \{f\in L^2(\mathbb S^1):f'\in L^2(\mathbb S^1) \}$. Thus, for $x\in[0,2\pi]$, we may write $$g(e^{\irm  x})=\sum_{m\in\mathbb Z}a_{m}e^{\irm mx}$$
for $a_m\in\mathbb C$ and we have $\sum_{m\in\mathbb Z}m^2|a_m|^2<\infty$ by Parseval's identity. It follows that $\sum_{m\in\mathbb Z}|a_m|<\infty$. Therefore, for $z\in\mathbb{D}(0,|\chi|^{-1})$, we obtain
\begin{align*}
    \sum_{n\geq0}|\chi|^ng(e^{\irm n\theta})z^n&=\sum_{m\in\mathbb Z}a_m\sum_{n\geq0}|\chi|^ne^{\irm mn\theta}z^n=\sum_{m\in\mathbb Z}a_m\frac{1}{1-e^
    {\irm m\theta}|\chi|z}.
\end{align*} 

Now we compute $a_m$. Since $g$ is piecewise linear, we can write 
$g(z)=b_jz+c_j\overline{z}$ when $\arg z\in[\alpha_j,\alpha_{j+1})$, where $0=\alpha_0<\alpha_1<\ldots<\alpha_{t}=2\pi$ are the break points on the unit circle $\mathbb S^1$ for some $t\in\mathbb N$. For $m\neq\pm1$, we have
\begin{align*}
  2\pi a_{m}&=\int_{0}^{2\pi}g(e^{\irm x})e^{-\irm m x}dx=\sum_{j=0}^{t-1}\int_{\alpha_j}^{\alpha_{j+1}}b_je^{\irm  (-m+1)x}+c_je^{\irm  (-m-1)x}dx\\&=\sum_{j=0}^{t-1}b_j\frac{e^{\irm (-m+1)\alpha_{j+1}}-e^{\irm (-m+1)\alpha_{j}}}{\irm(-m+1)}+\sum_{j=0}^{t-1}c_j\frac{e^{\irm (-m-1)\alpha_{j+1}}-e^{\irm (-m-1)\alpha_{j}}}{\irm(-m-1)}.
\end{align*}
For each $j=0,\ldots,t$, we denote $A_j=e^{\irm \alpha_j}$, then\begin{align*}
    a_{m}&=\frac{\irm(m+1)\sum_{j=0}^{t-1}b_j(A_{j+1}^{-m+1}-A_j^{-m+1})+\irm(m-1)\sum_{j=0}^{t-1}c_j(A_{j+1}^{-m-1}-A_j^{-m-1})}{2\pi(m^2-1)}\\
    &=\frac{d_m}{m^2-1}
\end{align*}
where $d_m$ is a power sum in $m$; thus, the sequence $\{d_m:m\in\mathbb Z\}$ satisfies a linear recurrence relation as wanted.
\end{proof}
\begin{remark}\label{degree of a power sum}
	It is known that every linear recurrence sequence can be represented as values of a power sum $\sum_{j=1}^h\rho_j^m P_j(m)$ for some distinct non-zero complex numbers $\rho_j$ and some non-zero polynomials $P_j$. Such a representation is unique. The \emph{degree} of such a power sum is defined as the maximal degree among the $P_j$.
	
	In Proposition~\ref{Fourier coefficients of continuous piecewise linear functions}, the sequence $d_m$ is represented by a power sum of degree at most 1.

\end{remark}

\section{Proof of Theorem~\ref{natural boundary of k-degree sequence} and its corollaries}\label{Section: natural boundary of generating series}

We have seen how to derive Corollary~\ref{the example of BDJ} from Theorem~\ref{natural boundary of k-degree sequence} in the introduction. In this section, we prove Theorem~\ref{natural boundary of k-degree sequence} and, subsequently,  Corollary~\ref{algebraically independence of generating functions}.

\subsection{Proof of Theorem~\ref{natural boundary of k-degree sequence}}

The proof of Theorem~\ref{natural boundary of k-degree sequence} is divided into three steps. 
\begin{itemize}
    \item Step 1: Expressing the $k$-degree sequence as values at $\Lambda^dA^n$ of a piecewise linear function $h\colon\mathrm{End}(\Lambda^kM_{\mathbb R})\rightarrow\mathbb R$. 
    \item Step 2: We cut the series into two parts: the dominant part $\Delta_{\Pi}(z)$ whose radius of convergence is $\lambda_k^{-1}$, and the remaining part $H(z)$ whose convergence radius is strictly larger than $\lambda_k^{-1}$. We show that $\Delta_{\Pi}(z)$ is not rational.
    \item Step 3: We use Proposition~\ref{Fourier coefficients of continuous piecewise linear functions} to express $\Delta_{\Pi}(z)$ as wanted. We then apply a result on the existence of the natural boundary of a power series to conclude, see Proposition~\ref{A general result on natural boundary}.
    \end{itemize}
\medskip
\noindent
\emph{Step 1}: Our computations, following~\cite{Lin-Wulcan-2014}, rely on the formulas for computing mixed volumes of polytopes that were developed in~\cite{Huber-Sturmfels-1995}.

Recall that the $k$-th exterior power $\Lambda^k M_{\mathbb R}$ is a real vector space of dimension $\binom{d}{k}$.

\begin{proposition}\label{a formula for the k-degree}
    For any polytope $P$ in $M_{\mathbb R}$ and for any $k \in \{1, \ldots , d\}$, there exists a piecewise linear function $h_P \colon \mathrm{End}(\Lambda^k M_{\mathbb R}) \rightarrow \mathbb R$, depending on $P$, such that the following holds.

For any endomorphism $u\colon M_{\mathbb R} \rightarrow M_{\mathbb R}$, one has
$$\mathrm{Vol} (u(P)[k], P[d-k]) = h_P(\Lambda^k u).$$
\end{proposition}
\begin{proof}
    We need some notation.
    
    Let $\mathcal{P}=(P_1,P_2)$ be a pair of polytopes in $M_\mathbb R$ such that $P_1+P_2$ has dimension $d$. A \emph{cell} of $\mathcal{P}$ is a pair $\mathcal{C}=(C_1,C_2)$ of polytopes $C_i\subset P_i$, $i=1,2$. We denote by $\#C_i$ the number of vertices of $C_i$, and by $C$ the Minkowski sum $C_1+C_2$. A \emph{fine mixed subdivision} of $\mathcal{P}$ is a collection of cells $\mathcal{S}=\left\{\mathcal{C}^{(1)},\ldots,\mathcal{C}^{(r)}\right\}$ such that $\dim C^{(j)}=d$, $\dim C_1^{(j)}+\dim C_{2}^{(j)}=d$, $\#(C_1^{(j)})+\#(C_2^{(j)})-2=d$ for all $1\leq j\leq r$, $C^{(j)}\cap C^{(j')}$ is a face of $C^{(j)}$ and $C^{(j')}$ for all $j,j'$, $\bigcup_{j}C^{(j)}=P_1+P_2$. 
    
     If $\mathcal{S}$ is a fine mixed subdivision of $\mathcal{P}$, then~\cite[Theorem~2.4]{Huber-Sturmfels-1995} says that 
\begin{equation}\label{a formula for mixed volume}\vol(P_1[k],P_2[d-k])=k!(d-k)!\sum_{\mathcal{C}^{(j)}\in\mathcal{S},\dim C_1^{(j)}=k,\dim C_2^{(j)}=d-k}\vol(C^{(j)}).
\end{equation}

Now, assume that $P$ has vertices $v_1,\ldots,v_s$ in $M_{\mathbb R}$ with $s=\#P$, then $u(P)$ has vertices $u(v_1),\ldots,u(v_s)$. We can apply~\eqref{a formula for mixed volume} to the fine mixed subdivision $\mathcal{S}$ of $(u(P),P)$ with cells of the form $$\mathcal{C}_{\mathcal{I}\mathcal{J}}\coloneqq (\mathrm{conv}(u(v_{i_0}),\ldots,u(v_{i_k})),\mathrm{conv}(v_{j_0},\ldots,v_{j_{d-k}}))$$ for some subsets of indices $\mathcal{I}=\{i_0,\ldots,i_k\}$ and $\mathcal{J}=\{j_0,\ldots,j_{d-k}\}$ in $\{1,\ldots,s\}$. Here, $\mathrm{conv}(\ldots)$ denotes the convex hull. Let $\mathcal{S}_k$ be the set of all such cells. It follows from~\eqref{a formula for mixed volume} that
$$\vol(u(P)[k],P[d-k])=k!(d-k)!\sum_{\mathcal{C}_{\mathcal{I}\mathcal{J}}\in \mathcal{S}_k}\vol(C_{\mathcal{I}\mathcal{J}}).$$

Since $C_{\mathcal{I}\mathcal{J}}=\mathrm{conv}(u(v_{i_0}),\ldots,u(v_{i_k}))+\mathrm{conv}(v_{j_0},\ldots,v_{j_{d-k}})$, we have $k!(d-k)!\vol(C_{\mathcal{I}\mathcal{J}})=|\det (D_{\mathcal{I}\mathcal{J}})|$ where $D_{\mathcal{I}\mathcal{J}}$ is the matrix with vectors $u(v_{i_1}-v_{i_0}),\ldots,u(v_{i_k}-v_{i_0}),$  $v_{j_1}-v_{j_0},\ldots,v_{j_{d-k}}-v_{j_0}$ as columns, see~\cite[Lemma~2.5]{Huber-Sturmfels-1995}. Observe that $\det (D_{\mathcal I\mathcal J})$ can be written as a $\mathbb Z$-linear combination of $k\times k$-minors of $u$, hence a $\mathbb Z$-linear combination of entries of $\Lambda^ku$. Let $\mathcal{C}_{\mathcal{I}\mathcal{J}}$ run over $\mathcal{S}_k$, the equations $\det (D_{\mathcal{I}\mathcal{J}})=0$ thus give us hyperplanes in $\mathrm{End}(\Lambda^k M_\mathbb R)$, hence decompose $\mathrm{End}(\Lambda^k M_\mathbb R)$ into a union of $\binom{d}{k}$-dimensional cones. Further, for any $u$ so that $\Lambda^ku$ lies on one of those cones (depending on, for each $\mathcal{C}_{\mathcal{I}\mathcal{J}}\in \mathcal{S}_k$, whether $\det(D_{\mathcal{I}\mathcal{J}})$ is non-negative or non-positive), the expression 
$$\mathrm{Vol}(u(P)[k],P[d-k])=\sum_{\mathcal{C}_{\mathcal{I}\mathcal{J}}\in \mathcal{S}_k}|\det(D_{\mathcal{\mathcal{I}\mathcal{J}}})|$$
gives us a linear function in terms of $\Lambda^ku$. Therefore, we obtain a piecewise linear function $h_P\colon\mathrm{End}(\Lambda^kM_{\mathbb R})\rightarrow\mathbb R$ with respect to some fan in $\mathrm{End}(\Lambda^k M_\mathbb R)$ such that  $\vol(u(P)[k],P[d-k])=h_P(\Lambda^ku)$ as desired.
\end{proof}

\medskip
\noindent
\emph{Step 2}: We may assume that $D$ is $\mathbb G_m^d$-invariant, in which case $D$ gives rise to a lattice polytope $P_
D$ in $M_\mathbb R$. By Proposition~\ref{a formula for the k-degree}, there is a piecewise linear function $h=d!h_{P_D}\colon \mathrm{End}(\Lambda^kM_\mathbb R)\rightarrow\mathbb R$ such that 
$$\mathrm{deg}_{D,k}(\varphi^n)=d!\vol((A^nP_D)[k],P_D[d-k])=h(\Lambda^kA^n)\text{ for all }n\geq0.$$

Recall that $1\leq k\leq d-1$ and the eigenvalues of $A\colon M_\mathbb R\rightarrow M_\mathbb R$ satisfy $|\mu_{k-1}|>|\mu_k|=|\mu_{k+1}|>|\mu_{k+2}|$  and $\mu_{k},\mu_{k+1}$ are complex conjugates whose arguments are incommensurable with $2\pi$. Therefore, $\lambda_{k-1}=|\mu_1\mu_2\cdots\mu_{k-1}|$, $\lambda_k=|\mu_1\mu_2\cdots\mu_{k}|$ by~\cite[Corollary~B]{Favre-Wulcan-2012}, and the matrix $\Lambda^kA$ has two dominant eigenvalues $\chi_{1}\coloneqq \mu_1\cdots\mu_{k-1}\mu_{k}$ and $\chi_{2}\coloneqq \mu_1\cdots\mu_{k-1}\mu_{k+1}$. The other eigenvalues of $\Lambda^kA$ have moduli less than or equal to $\lambda\coloneqq \lambda_{k-1}|\mu_{k+2}|<\lambda_k$.  Observe that the product $\mu_1\mu_2\cdots\mu_{k-1}$ is real; thus, $\chi_1$ and $\chi_2$ are complex conjugates
whose arguments, say $\theta$ and $-\theta$ respectively, are incommensurable with $2\pi$. We obtain:

\begin{lemma}
   There exists a $\Lambda^k A$-invariant splitting $\Lambda^k M_\mathbb R = V \oplus W$
with $\mathrm{dim} (V) = 2$, and
    \begin{itemize}
        \item the eigenvalues of the restriction of $\Lambda^kA$ to $V$ are $\chi_1$ and $\chi_{2}$,
        \item the spectral norm of the restriction on $W$ is $\lambda$.
    \end{itemize}
\end{lemma}  

For $n\geq0$, we write $\Lambda^kA^n|_V$ for the restriction of $\Lambda^kA^n$ to $V$. Observe that the subspace $\Pi$ of endomorphisms of $V$ that commute with $\Lambda^kA|_V$ is equal to $\mathbb R \mathrm{Id} \oplus \mathbb R \Lambda^kA|_V$, hence is two-dimensional. We remark that $\Lambda^kA^n|_V\in\Pi$ for all $n\geq0$.

We can identify $\Pi$ with its image in $\mathrm{End}(\Lambda^kM_\mathbb R)$ under the map
$u \mapsto u \oplus 0_W$. To simplify notation, we denote by $h\colon \Pi \rightarrow \mathbb R$
the piecewise linear function obtained by restricting $h$ to $\Pi$, so that
$h(\Lambda^kA^n|_V)$ is the value of $h$ at the endomorphism $\Lambda^kA^n|_V \oplus 0_W$.

Since $h$ is Lipschitz by Lemma~\ref{piecewise linear functions are Lipschitz}, we have 
\[\deg_{D,k}(\varphi^n)=h(\Lambda^kA^n)=h(\Lambda^kA^n|_V)+O(\lambda^n);\] hence $h(\Lambda^kA^n|_V)\asymp \lambda_k^n$ since $\deg_{D,k}(\varphi^n)\asymp \lambda_k^n$ by~\cite[Theorem A]{Favre-Wulcan-2012}. It follows that
$$\Delta_{k,\varphi,D}(z)=\Delta_{\Pi}(z)+H(z)$$
where $\Delta_{\Pi}(z)=\sum_{n\geq0}h(\Lambda^kA^n|_V)z^n$ converges inside $\mathbb D(0,\lambda_k^{-1})$ and  $H(z)$ is a power series whose convergence radius is greater than or equal to $\lambda^{-1}>\lambda_k^{-1}$. 

\begin{lemma}\label{k-degree sequence does not satisfy linear recurrence}
The sequence $h(\Lambda^kA^n|_V), n\geq0$, does not satisfy any linear recurrence relation.   
\end{lemma}
\begin{proof}
    Assume the contrary. Since $h\colon\Pi\rightarrow\mathbb R$ is piecewise linear with respect to some fan, we can find a cone $\sigma$ in $\Pi$ so that $\Lambda ^kA^n|_V$ lies in $\sigma$ for infinitely many $n$. In other words, there is some linear function $h_{\sigma}\colon\Pi\rightarrow\mathbb R$ so that $h(\Lambda^kA^n|_V)=h_\sigma(\Lambda^kA^n|_V)$ for infinitely many $n$. We denote $\beta_n=h_\sigma(\Lambda^kA^n|_V)$ for $n\geq0$. 
    
    The Cayley-Hamilton theorem implies that the sequence of matrices $\Lambda^kA^n|_V$, $n\geq0$, satisfies the linear recurrence relation given by the characteristic polynomial of $\Lambda^kA|_V$. Therefore, the sequence $\beta_n$ satisfies the same linear recurrence relation. Since the sequence $h(\Lambda^kA^n|_V)$ also satisfies a linear recurrence relation and coincides with the sequence $\beta_n$ for infinitely many $n$, these two sequences must coincide on some arithmetic progression $a+b\mathbb N$ with $b>0$ by the Skolem--Mahler--Lech  theorem, see, e.g.,~\cite[Chapter~4 Exercise~4]{Stanley-Enumerative-Combinatorics-2011}. Since $h(\Lambda^kA^n|_V)\asymp \lambda_k^n$, we have $\beta_n\asymp\lambda_k^n$ for $n\in a+b\mathbb N$. In particular, $\beta_n>0$ for $n\in a+b\mathbb N$ sufficiently large. 

Recall that $\theta=\arg\chi_1\not\in2\pi\mathbb Q$ and $b\neq0$. We may write 
$$\beta_{a+bn}=h_\sigma(\Lambda^kA^{a+bn}|_V)=C\mathrm{Re}(\chi_1^{a+bn})+D\mathrm{Im}(\chi_1^{a+bn})=\lambda_k^{a+bn}\mathrm{Re}((C-\irm D)e^{\irm (a+bn)\theta})$$ 
 for some $C,D\in\mathbb R$, not both equal to zero. Since $\{e^{\irm (a+bn)\theta}:n\in\mathbb N\}$ is dense in the unit circle $\mathbb S^1$, there are infinitely many $n\in\mathbb N$ such that  $\mathrm{Re}((C-\irm D)e^{\irm (a+bn)\theta})<0$, which is absurd.
\end{proof}
 
It follows from~\cite[Theorem~4.1.1]{Stanley-Enumerative-Combinatorics-2011} that $\Delta_{\Pi}(z)$ is not a rational function.

\medskip
\noindent
\emph{Step 3}: Recall that $\theta$ is the argument of $\chi_1=\mu_1\cdots\mu_{k-1}\mu_{k}$ and $\lambda=|\mu_1\cdots\mu_{k-1}\mu_{k+2}|<\lambda_k$.
\begin{lemma}\label{an expression of Delta}
Under the assumption of Theorem~\ref{natural boundary of k-degree sequence}, there are complex numbers $\{a_m\}_{m\in\mathbb Z}$ such that $\sum_{m\in\mathbb Z}|a_{m}|<\infty$ and
     $$\Delta_{\Pi}(z)=\sum_{m\in\mathbb Z}\frac{a_{m}}{1-e^{\irm m\theta}\lambda_kz}$$ for all  $z\in\mathbb D(0,\lambda_k^{-1})$. Moreover, there exists a linear recurrence sequence $\{d_{m}:m\in\mathbb Z\}$ such that $a_{m}=\frac{d_{m}}{m^2-1}$ for all $m\neq\pm1$. In addition, there exists an arithmetic progression $a+b\mathbb N$ with $b\neq0$ such that $d_m\neq0$ for all $m \in a+b\mathbb N$.  
\end{lemma}
\begin{proof}
We endow the plane $V$ with a conformal structure so that
it gets identified with the complex plane $\mathbb C$, and $\Lambda^k A|_V$ becomes the multiplication
by $\chi_1$. Via such identifications, the piecewise linear function $h:\Pi\rightarrow \mathbb R$ can be viewed as a piecewise linear function $h\colon\mathbb R^2\rightarrow\mathbb R$ for which $h(\Lambda^k A^n|_V)=h(\mathrm{Re}(\chi_1^n),\mathrm{Im}(\chi_1^n))$ for all $n\geq0$. By Proposition~\ref{Fourier coefficients of continuous piecewise linear functions}, we can write
$$\Delta_{\Pi}(z)=\sum_{m\in\mathbb Z}\frac{a_m}{1-e^
    {\irm m\theta}\lambda_kz}$$ for $z\in\mathbb D(0,\lambda_k^{-1})$. In addition, for $m\neq\pm1$ we have  $a_m=\frac{d_m}{m^2-1}$  for some linear recurrence sequence $\{d_m:m\in\mathbb Z\}$ as desired.   

    Since $\Delta_\Pi(z)$ is not rational, there are infinitely many $m\in\mathbb Z$ so that $d_m\neq0$. The last assertion then follows from the Skolem--Mahler--Lech  theorem~\cite[Chapter~4 Exercise 4]{Stanley-Enumerative-Combinatorics-2011}.
\end{proof}
Since $H(z)$ converges inside a disk properly containing the disk $\mathbb D(0,\lambda_k^{-1})$, it remains to show that $\Delta_\Pi(z)$ has 
$\mathcal{C}(0,\lambda_k^{-1})$ as a natural boundary. Since $\{e^{\irm (a+bn)\theta}\lambda_k^{-1}:n\in\mathbb N\}$ is dense in $\mathcal{C}(0,\lambda_k^{-1})$, we complete the proof of Theorem~\ref{natural boundary of k-degree sequence} by applying the following result.

\begin{lemma}\label{A general result on natural boundary}
    Let $\{p_m:m\in \mathbb{Z}\}\subset\mathcal C(0,\lambda_k)$. Let  $\{b_{m}\}_{m\in\mathbb Z}$ be complex numbers such that $\sum_{m\in \mathbb Z}|b_{m}|<\infty$. Consider the power series $\Delta(z)\coloneqq \sum_{m\in \mathbb{Z}}\frac{b_{m}}{1-p_mz}.$  Then we have  $$\lim_{\rho\to1^{-}}(1-\rho)\Delta(p_m^{-1}\rho)=b_{m}\text{ }\text{ for all }m\in\mathbb Z.$$
    If we assume further that the set of points $p_m\in\mathcal{C}(0,\lambda_k)$ satisfying $b_m\neq0$ is dense in $\mathcal{C}(0,\lambda_k)$, then $\Delta(z)$ has $\mathcal{C}(0,\lambda_k^{-1})$ as a natural boundary.
\end{lemma}

\begin{proof}[Proof of Lemma~\ref{A general result on natural boundary}]
   The hypothesis yields that $\Delta(z)$ converges in $\mathbb D(0,\lambda_k^{-1})$. Now, we fix some $m_0\in\mathbb Z$, and take $z= p_{m_0}^{-1}\rho$ with $\rho\in\mathbb R$, $0<\rho<1$. Then 
$$\Delta(p_{m_0}^{-1}\rho)=\sum_{m\in\mathbb Z\setminus\{m_0\}}\frac{b_m}{1-p_mp_{m_0}^{-1}\rho}+\frac{b_{m_0}}{1-\rho}.$$
Multiplying both sides by $(1-\rho)$ yields
$$(1-\rho)\Delta(p_{m_0}^{-1}\rho)=\sum_{m\in\mathbb Z\setminus\{m_0\}}{b_m}\frac{1-\rho}{1-p_mp_{m_0}^{-1}\rho}+b_{m_0}.$$
For every $m\neq m_0$, by triangle inequality we have $$\Bigg|\frac{1-\rho}{1-p_mp_{m_0}^{-1}\rho}\Bigg|\leq \frac{1-\rho}{1-|p_mp_{m_0}^{-1}\rho|}\leq\frac{1-\rho}{1-\rho}\leq 1\text{ }\text{ for }0<\rho<1.$$
Thus, it follows from $\sum_{m\in\mathbb Z}|b_m|<\infty$ that 
$$\sum_{m\in\mathbb Z\setminus\{m_0\}}{|b_m|}\Bigg|\frac{1-\rho}{1-p_mp_{m_0}^{-1}\rho}\Bigg|\leq \sum_{m\in\mathbb Z\setminus\{m_0\}}|b_m|<\infty.$$
The Lebesgue dominated convergence theorem implies that \begin{align*}
    \lim_{\rho\to 1^-}\sum_{m\in\mathbb Z\setminus\{m_0\}}{b_m}\frac{1-\rho}{1-p_mp_{m_0}^{-1}\rho}=\sum_{m\in\mathbb Z\setminus\{m_0\}}\lim_{\rho\to1^-}{b_m}\frac{1-\rho}{1-p_mp_{m_0}^{-1}\rho}=0.
\end{align*}
Therefore 
$\lim_{\rho\to 1^-}(1-\rho)\Delta(p_{m_0}^{-1}\rho)=b_{m_0}.$

In addition, if $b_{m_0}\neq0$, then $\lim_{\rho\to 1^-}\Delta(p_{m_0}^{-1}\rho)=\infty$. Thus, if the set of points $p_m$ satisfying $b_m\neq0$ is dense in $\mathcal{C}(0,\lambda_k)$, then $\Delta(z)$ cannot be
analytically continued at any point on $\mathcal{C}(0,{\lambda_k}^{-1})$ as wanted.

\end{proof}

\subsection{Proof of Corollary~\ref{algebraically independence of generating functions}}

Recall that we want to prove $\mathrm{tr.deg}_{\mathbb C(z)}(\Delta_{k,f,D}(z):(X,K,D,f)\in \mathfrak M)=\infty$ for all $k\geq1$. We fix $k\geq1$. Let $X(\Sigma)$ be any projective toric variety of dimension $k+1$ over some algebraically closed field, and $D$ an ample Cartier divisor on $X(\Sigma)$. For instance, we can take $\mathbb P^{k+1}$ with $\mathcal{O}(1)$. It suffices to construct infinitely many dominant monomial self-maps of $X(\Sigma)$ whose associated generating series of their $k$-degree sequences are algebraically independent over $\mathbb C(z)$. 

\begin{lemma}\label{Lemma: algebraically independent of delta}
     Let $\varphi_1,\ldots,\varphi_s$ be dominant monomial self-maps of $X(\Sigma)$. Suppose that all $\varphi_i$ satisfy the conditions in Theorem~\ref{natural boundary of k-degree sequence}. Suppose further that their $k$-dynamical degrees are pairwise distinct. Then the series $\Delta_{k,\varphi_1,D},\ldots,\Delta_{k,\varphi_s,D}$ are algebraically independent over $\mathbb C(z)$.  
\end{lemma}


Granted this lemma, we can construct our desired monomial maps as follows. First, we choose any monomial self-map $\varphi$ satisfying the assumptions of Theorem~\ref{natural boundary of k-degree sequence}. It follows that $\lambda_k(\varphi)>1$. Thus the $k$-dynamical degrees $\lambda_k(\varphi^\ell)=\lambda_k(\varphi)^\ell$, $\ell\geq1$, are pairwise distinct. Lemma~\ref{Lemma: algebraically independent of delta} yields that 
$\mathrm{tr.deg}_{\mathbb{C}(z)} \bigl(\Delta_{k,\varphi^\ell,D}(z) \bigm| \ell \geq 1\bigr) = \infty,$
as desired.

It remains to prove Lemma~\ref{Lemma: algebraically independent of delta}. By Theorem~\ref{natural boundary of k-degree sequence}, the series $\Delta_{k,\varphi_1,D_1},\ldots,\Delta_{k,\varphi_s,D_s}$ have circles of convergence as natural boundaries. The lemma then follows from the following observation.

\begin{lemma}
    Let $\Delta_1,\ldots,\Delta_s$ be power series in $\mathbb C[[z]]$ that have circles of pairwise distinct positive radii as natural boundaries. Then they are algebraically independent over $\mathbb C(z)$.
\end{lemma}
\begin{proof}
    Let $\rho_1,\ldots,\rho_s$ be their radii of convergence, ordered so that $0<\rho_1<\rho_2<\ldots<\rho_s$. Suppose by contradiction that there is a polynomial relation between these functions. We may suppose that this relation is minimal in terms of the number of power series $\Delta_i$ that appear. We see that
    there are polynomials $P_0,\ldots,P_r\in\mathbb C[z][X_2,\ldots,X_\ell]$ for some $1\leq \ell\leq s$ and non-zero $P_r$ such that $P_r(\Delta_2,\ldots,\Delta_\ell)\neq0$ and
    \begin{equation}\label{algebraic relation between Delta}
\Delta_1^rP_r(\Delta_2,\ldots,\Delta_\ell)+\cdots+P_0(\Delta_2,\ldots,\Delta_\ell)=0.
    \end{equation}
    The left-hand side of~\eqref{algebraic relation between Delta} can be viewed as a polynomial in $\Delta_1$ with coefficients in the ring of holomorphic functions inside the domain $\mathbb D(0,\rho_2)$. We denote by $\Delta$ its discriminant when considering $\Delta_1$ as a variable. It is a holomorphic function inside $\mathbb D(0,\rho_2)$.
    
    Recall that $\Delta_1$ has the circle $\mathcal{C}(0,\rho_1)$ as a natural boundary. By the identity theorem, we can find a point $t\in\mathcal{C}(0,\rho_1)$ such that $P_r(\Delta_2(t),\ldots,\Delta_\ell(t))\neq0$ and $\Delta(t)\neq0$. Then $P_r(\Delta_2(z),\ldots,\Delta_\ell(z))\neq0$ and $\Delta(z)\neq0$ around some disk $\mathbb D$ centered at $t$. Now, there exists exactly $r$ complex numbers $c$ satisfying the algebraic equation $$c^rP_r(\Delta_2(t),\ldots,\Delta_\ell(t))+\cdots+P_0(\Delta_2(t),\ldots,\Delta_\ell(t))=0.$$
By the Cauchy--Kovalevskaya theorem~\cite[Theorem~9.4.8]{Hormander-book-1983}, for each $c$, there is an open disk $\mathbb D_0\subseteq\mathbb D$ centered at $t$ and a holomorphic function $y$ on $\mathbb D_0$ such that $y(t)=c$. By the uniqueness of solutions of the previous initial value problem, we deduce that $\Delta_1$ must equal some $y$ on $\mathbb D_0\cap \mathbb D(0,\rho_1)$. In other words, there exists an analytic extension of $\Delta_1$ to $\mathbb D_0\cap \mathbb D(0,\rho_1)$, which is absurd. 
\end{proof}

 \section{The surface case}\label{The surface case.}

In this section, we restrict our attention to monomial surface maps. The advantage is that the eigenvalues of the associated $2\times 2$ matrices are easier to control than in the general case.

\subsection{Proof of the dichotomy part of Theorem~\ref{dichotomy of the generating power series for toric surfaces}}\label{subsection on dichotomy}

There are two possibilities for the eigenvalues of $A\in M_2(\mathbb Z)$ of non-zero determinant. 
\begin{enumerate}
    
    \item If $\overline{\mu_2}=\mu_1$ and $\mu_1/\mu_2$ is not a root of unity, then it follows from Theorem~\ref{natural boundary of k-degree sequence} that $\Delta_{1,\varphi,D}$ has the desired form. In addition, $\Delta_{1,\varphi,D}$ has the circle $\mathcal{C}(0,\lambda_1^{-1})$ as a natural boundary.
    \item In the remaining case,  by~\cite[Th\'eor\`eme principal]{Favre-2003}, $\varphi$ can be lifted to a $1$-stable monomial map. Here, we note that toric surfaces are always simplicial. Therefore, $\Delta_{1,\varphi,D}$ is rational thanks to~\cite[Corollary~2.2]{Diller-Favre-2001}\footnote{Although the article was written for surfaces over $\mathbb C$, the result works with a proof unchanged over any algebraically closed field.}.   
\end{enumerate}

   This finishes the proof of the dichotomy.
\begin{remark}\label{a closed form for dm}
   In the first case, we have $$\Delta_{1,\varphi,D}(z)=\sum_{m\in\mathbb Z}\frac{a_m}{1-e^{\irm m\theta}\lambda_1z}$$ where $\theta=\arg\mu_1$ and $a_m=\frac{d_m}{m^2-1}$ for some linear recurrence sequence $\{d_m:m\in\mathbb Z\}$ with $\sum_{m\in\mathbb Z}|a_m|<\infty$. Further, the linear recurrence sequence $d_m$ is represented by a power sum of degree at most 1, see Remark~\ref{degree of a power sum}. This observation will be useful in the proof of Theorem~\ref{when two degree sequences coincide}.
\end{remark}

\begin{remark}
    We propose a direct proof of Theorem~\ref{dichotomy of the generating power series for toric surfaces} without relying on~\cite{Favre-2003}. When $\mu_1>\mu_2>0$, since $d=2$, we have
    $$\Delta_{1,\varphi,D}(z)=\Delta_\Pi(z)=\sum_{n\geq0}h(\mu_1^n,\mu_2^n)=\sum_{n\geq0}\mu_1^nh\Big(1,\Big(\frac{\mu_2}{\mu_1}\Big)^n\Big)z^n$$ for some piecewise linear function $h\colon \mathbb R^2\rightarrow\mathbb R$ with respect to some (not necessarily rational) fan. Since $(\mu_{2}/\mu_1)^n$ tends to 0 as $n$ goes to infinity, there is some $N_0\in\mathbb N$ so that the points $(1,({\mu_2}/{\mu_1})^n)$ belong to some two-dimensional cone of the fan for all $n\geq N_0$. Thus, 
    there are some $p,q\in\mathbb R$ such that for $n\geq N_0$, we have  $h(1, (\mu_2/\mu_1)^n)=p+q(\mu_2/\mu_1)^n.$
    Therefore  $$\Delta_{1,\varphi,D}(z)=\sum_{n=0}^{N_0}h(\mu_1^n,\mu_2^n)z^n+\sum_{n\geq N_0}(p\mu_1^n+q\mu_2^n)z^n\in\mathbb C(z).$$
    Other cases for the eigenvalues $\mu_1,\mu_2$ can be treated similarly.
\end{remark}

\begin{remark}
In light of Remark~\ref{remark1}, we can ask whether, when $\overline{\mu_2}=\mu_1$ with $\mu_2/\mu_1$ not being a root of unity,  $\Delta_{1,\varphi,D}$ is differentially transcendental over $\mathbb C(z)$. Although we are not able to prove it for $\Delta_{1,\varphi,D}$, we can prove it for $$\Delta_0(z)\coloneqq \sum_{m\in\mathbb Z}\frac{(m^2-1)a_m}{1-e^{\irm m\theta}\lambda_1z}=\frac{-a_{0}}{1-\lambda_1z}+\sum_{|m|>1}\frac{d_m}{1-e^{\irm m\theta}\lambda_1z}.$$
Indeed, since both $\{d_m\}_{m>1}$ and $\{d_m\}_{m<-1}$ are linear recurrence sequences, we deduce that the series $\Delta_0(z)$ satisfies a non-trivial linear $e^{\irm\theta}$-difference equation with coefficients in $\mathbb C(z)$. Since $e^{\irm\theta}$ is not a root of unity, it follows from~\cite[Theorem~1.2]{Adamczewski-Dreyfus-Hardouin-2021} that either $\Delta_0\in \bigcup_{j\geq1}\mathbb C(z^{1/j})$, or it is differentially transcendental. Similar to $\Delta_{1,\varphi,D}(z)$, we observe that $\Delta_0(z)$ also admits the circle $\mathcal{C}(0,\lambda_1^{-1})$ as a natural boundary; in particular, $\Delta_0(z)\not\in \bigcup_{j\geq1}\mathbb C(z^{1/j})$. Therefore $\Delta_0(z)$ is differentially transcendental as wanted. We expect that there is a functional equation between $\Delta_{1,\varphi,D}$ and $\Delta_0$ from which one can deduce the differential transcendence for $\Delta_{1,\varphi,D}$. 
\end{remark}

\subsection{Proof of the reduction part of Theorem~\ref{dichotomy of the generating power series for toric surfaces}}
\noindent

Our goal is to prove that when the eigenvalues $\mu_1$ and $\mu_2$ of the matrix $A$ are complex conjugates and $\mu_1/\mu_2$ is not a root of unity, then for all prime numbers $p$ sufficiently large, the reduction $\Delta_{1,\varphi,D}(z)\bmod p$ is transcendental over $\mathbb F_p(z)$. We need the following result from automata theory. 
\begin{lemma}\label{algebraic implies weakly periodic}
     Let $\mathbb F_q$ be a finite field of characteristic $p>0$. If $\sum_{n\geq0}s_nz^n\in \mathbb F_q((z))$ is algebraic over $\mathbb F_q(z)$, then the sequence $(s_n)_{n\geq0}$ is weakly periodic in the sense that for any integers $a>0$ and $b\geq0$, there exist $k>0$ and $r>r'\geq0$ such that
$s_{a(kn+r)+b}=s_{a(kn+r')+b}\text{ for all }n\geq0.$
\end{lemma}
\begin{proof}
   This follows from Christol's theorem~\cite[Th\'eor\`eme~1]{Christol-1979}\footnote{Although the article stated the result only for $\mathbb{F}_p$, it holds verbatim for $\mathbb{F}_q$, cf.~\cite[Th\'eor\`eme~1]{Christol-Kamae-France-Rauzy-1980}.} and~\cite[Lemma~2.1]{Byszewski-Konieczny-2020}.
\end{proof}

Therefore, we need to show that for primes $p$ sufficiently large, the reduction modulo $p$ of the degree sequence is not weakly periodic. In the proof of Theorem~\ref{natural boundary of k-degree sequence}, we have interpreted the degree sequence as values of a piecewise linear function. However, we need an additional property for such a function, namely, its convexity. 

 We may assume that $D$ is an ample $\mathbb G_m^2$-invariant Cartier divisor. Recall that $P_D\subset M_\mathbb R$ is the polyhedron associated with $D$, and $h_D\colon N_\mathbb R\rightarrow \mathbb R$ is its associated strictly convex piecewise linear function with respect to the fan $\Sigma$. We note that when restricted to each cone of $\Sigma$, $h_D$ is a linear function with coefficients in $\mathbb Z$. We recall also that the subspace $\Pi$ of endomorphisms of $M_\mathbb R$ that commute with $A$ is equal to $\mathbb R\mathrm{Id}\oplus\mathbb RA$, which is two-dimensional. 
\begin{lemma}\label{a formula for the 1-degree}
 There exists a convex piecewise linear function $h\colon\Pi\rightarrow \mathbb R$ with respect to some complete fan in $\Pi$, which is not linear, such that for any endomorphism $u\in\Pi$, one has
 $2\mathrm{Vol}(u(P_D),P_D)=h(u).$ In particular,  we have $$\deg_{D,1}(\varphi^n)=h(A^n)\text{ for all }n\geq0.$$
 Further, for each cone $\sigma$ over which $h$ is linear, the coefficients of $h_\sigma$ are in $\mathbb Z$.  
\end{lemma}
\begin{proof}
    We will use the \emph{surface area measure} to compute the mixed volume, following~\cite{Schneider-book-2013}. We fix a scalar product for which the given basis of $N$ forms an orthogonal basis, which gives rise to a Riemannian metric on $\mathbb S^1\subset N_\mathbb R$. If $Q$ is a compact polyhedron in $M_\mathbb R$, the surface area measure $\mathcal{S}(Q)$ on $\mathbb S^1\subset N_\mathbb R$ is defined as $\sum_{v}\mathrm{Vol}(F_v)\delta_v$ where $v$ runs over the outer unit normal vectors of the facets $F_v$ of $Q$, and $\delta_v$ denotes the Dirac measure at $v$. Here, the volume of each facet is its one-dimensional volume, in other words, its length. By~\cite[Formulae~(5.19) and~(5.23)]{Schneider-book-2013}, we have
    $$\mathrm{Vol}(u(P_D),P_D)=\frac{1}{2}\int_{\mathbb S^1}(h_D \circ \check{u}) d\mathcal{S}(P_D)=\frac{1}{2}\sum_{v}\mathrm{Vol}(F_v)h_D(\check{u}(v))$$where $\check{u}\colon N_\mathbb R\rightarrow N_\mathbb R$ is the dual linear map of $u$, and $v$ runs over the outer unit normal vectors of the facets $F_v$ of $P_D$. We set $h(u)\coloneqq \sum_{v}\mathrm{Vol}(F_v)h_D(\check{u}(v))$ for $u\in \Pi$. 

Since all coefficients $\mathrm{Vol}(F_v)$ are positive, it is sufficient to show that for any $v$, the function $h_v(u)\coloneqq h_D(\check{u}(v))$ is piecewise linear and convex. Recall the fact that every convex piecewise linear function can be represented as the maximum of linear forms; thus, $h_D=\max\{l_1,\ldots,l_n\}$ for some linear forms $l_i$ with coefficients in $\mathbb Z$. It follows that $$h_v(u)=\max\{L_1(u),\ldots,L_n(u)\}$$ where $L_i(u)=l_i(\check{u}(v))$ is a linear form on the space $\Pi$. Thus $h_v$ is a convex piecewise linear function on $\Pi$ for every $v$. In addition, since the vectors $\mathrm{Vol}(F_v)v$ have coefficients in $\mathbb Z$ (as $P_D$ is a lattice polytope), the linear functions $\mathrm{Vol}(F_v)L_i(u)=l_i(\check{u}(\mathrm{Vol}(F_v)v))$ have coefficients in $\mathbb Z$. Therefore,
the piecewise linear functions $\mathrm{Vol}(F_v)h_v=\max\{\mathrm{Vol}(F_v)L_i\}$ also have coefficients in $\mathbb Z$, and so does $h$.

Since  $A^n\in\Pi$ for all $n\geq0$, we have
$\deg_{D,1}(\varphi^n)=2!\mathrm{Vol}(A^nP_D,P_D)=h(A^n)\text{ for all }n\geq0.$ Finally, if $h$ is linear, then $\Delta_{1,\varphi,D}(z)=\sum_{n\geq0}h(A^n)z^n$ would be rational which contradicts our result in \S\ref{subsection on dichotomy}.
\end{proof}
\begin{remark}
    This map $h$ is the same as the map $h$ constructed in Step 3 of the proof of Theorem~\ref{natural boundary of k-degree sequence} when $d=2$ and $k=1$, but now we have an additional convex property.
\end{remark}

Since $A$ has complex conjugates $\mu_1$ and $\mu_2$ as eigenvalues, we can find $Q\in\mathrm{GL}_2(\mathbb Q(\mu_1))$ such that $QAQ^{-1}=\begin{pmatrix}
    \mu_1&0\\0&\mu_2
\end{pmatrix}$. By making a linear change by $Q$, we can endow $M_\mathbb R$ with a conformal structure and $A$ becomes the multiplication by $\mu_1$. Via this identification, we can view $h\colon\Pi\rightarrow\mathbb R$ as a convex piecewise linear function $h\colon \mathbb C\rightarrow\mathbb R$, which is not linear, such that $\deg_{D,1}(\varphi^n)=h(\mu_1^n)$. We note that $h$ now has coefficients in $\mathbb Q(\mu_1)$ on each cone. Namely, there is an integer $s\geq1$ together with a decomposition of $[0,1]=\bigcup_{i=1}^sI_i$ into closed subsets, and there exist complex numbers $c_1,\ldots,c_s$ in $\mathbb Q(\mu_1)$ such that $h(z)=c_iz+\overline{c_i}\overline{z}$ whenever $\arg z\in 2\pi I_i$. Here, the closed sets $I_i$ are subsets of the form $[\alpha,\beta]$ or of the form  $[0,1]\setminus(\alpha,\beta)$ for some $0\leq \alpha<\beta\leq 1$ with $\beta-\alpha<1$. Via the identification $[0,1]\rightarrow\mathbb S^1,x\mapsto e^{\irm2\pi x}$, we may assume that the sets $I_i$, for $i=1,\ldots,s$, are organized so that they correspond to consecutive closed connected subsets covering the unit circle $\mathbb S^1$ in clockwise order. For each $i$, we introduce the cone $\sigma_i\coloneqq \{\lambda e^{\irm 2\pi x}:\lambda\geq0,x\in I_i\}$ in $\mathbb C$, then $h(z)=c_iz+\overline{c_i}\overline{z}$ whenever  $ z\in \sigma_i$.

\begin{lemma}\label{convex+linear implies strictly convex}
    If there are $1\leq i< j\leq s$ such that $c_i=c_j$, then either $c_{\ell}=c_i$ for all integers $\ell\in (i,j)$ or $c_{\ell}=c_i$ for all integers $\ell\not\in[i,j]$.
\end{lemma}
\begin{proof}
Replacing $h(z)$ by $h(z)-(c_iz+\overline{c_iz})$, we may assume that $c_i=c_j=0$. Thus $h(z)=0$ for all $z\in C_i\cup C_j$. 

We set $\sigma$ to be the convex hull of $\sigma_i$ and $\sigma_j$. Then either $\bigcup_{ \ell\in (i,j)}\sigma_{\ell}\subset \sigma$ or $\bigcup_{ \ell\not\in [i,j]}\sigma_{\ell}\subset \sigma$ (or both). We assume the former case and aim to prove that $c_{\ell}=c_i$ for all $\ell\in (i,j)$; the latter case is proved similarly.

Under our assumptions, we have $h(z)\leq 0$ for all $z\in \sigma$ thanks to the convexity of $h$. Let $z_i$ be any point in the interior of $\sigma_i$. For any $z\in \sigma$, we then can find some $0<t<1$ such that $tz_i+(1-t)z\in \sigma_i$, whence 
$$0=h(tz_i+(1-t)z)\leq th(z_i)+(1-t)h(z)=(1-t)h(z),$$
which implies that $ h(z)\geq0$. We deduce that $h(z)=0$ for all $z\in\sigma$. This follows that for any integer $\ell$ in $(i,j)$, we have $h(z)=0$ for $z\in \sigma_{\ell}$, hence $c_{\ell}=0$ as wanted.
\end{proof}

We can merge consecutive closed subsets $I_i$ on which $h$ has the same coefficients. Thanks to Lemma~\ref{convex+linear implies strictly convex}, we may then assume that the coefficients $c_1,\ldots,c_s$ are pairwise distinct. Since $h$ is not linear, we have $s>1$.
\begin{lemma}
    For every $1\leq i< j\leq s$, there is at most one $n\geq0$ such that $c_i\mu_1^n+\overline{c_i}\mu_2^n=c_j\mu_1^n+\overline{c_j}\mu_2^n$.
\end{lemma}
\begin{proof}
  For every $0\leq n<m$, we note that two vectors $(\mu_1^n,\mu_2^n)$ and $(\mu_1^m,\mu_2^m)$ are linearly independent over $\mathbb C$, otherwise $\mu_1^{m-n}=\mu_2^{m-n}$, which is absurd since $\mu_1/\mu_2$ is not a root of unity. Since the kernel of the non-trivial linear form $(c_i-c_j)z_1+(\overline{c_i}-\overline{c_j})z_2$ has dimension one over $\mathbb C$, the lemma follows.
  \end{proof}
 We conclude that there is some positive integer $b$ such that $c_{1}\mu_1^b+\overline{c_{1}}\mu_2^b\neq c_{j}\mu_1^b+\overline{c_{j}}\mu_2^b$ for all $1<j\leq s$. To simplify notation, let $K=\mathbb Q(\mu_1)$, consider the finite subset of non-zero elements in $K$
 $$\mathcal{F}=\{\mu_1,c_i\text{ for }1\leq i\leq s,\text{ and }c_{1}\mu_1^b+\overline{c_{1}}\mu_2^b- c_{i}\mu_1^b-\overline{c_{i}}\mu_2^b\text{ for }1<j\leq s\}.$$
  The set of places (ie., of equivalence classes of non-trivial norms) $|\cdot|_v$ on $K$ such that $|x|_v\neq1$ for some $x\in\mathcal{F}$ is finite. We infer that for all prime integers $p$ sufficiently large, and for any place $|\cdot|_v$ extending the $p$-adic norm to $K$, we have $|x|_v=1$ for all $x\in\mathcal{F}$. Fix any such prime $p$ and any place $|\cdot|_v$ extending the $p$-adic norm. We denote $K_v^\circ=\{x\in K:|x|_v\leq1\}$ and $K_v^{\circ\circ}=\{x\in K:|x|_v<1\}$, and observe that the residue field $\tilde{K}_v=K_v^\circ/K_v^{\circ\circ}$ is a finite extension of $\mathbb F_p$ with cardinality $q$. 

Assume by contradiction that $\Delta_{1,\varphi,D}$ modulo $p$ is algebraic over $\mathbb F_p(z)$.  Thanks to the natural embedding $\mathbb F_p\xhookrightarrow{}\tilde{K}_v$, the reduction $\Delta_{1,\varphi,D}(z)\bmod K_v^{\circ\circ}=\sum_{n\geq0}(h(\mu_1^n)\bmod K_v^{\circ\circ})z^n\in\tilde{K}_v[[z]]$ is also algebraic over $\mathbb F_p(z)$, hence over $\tilde{K}_v(z)$. It follows that the sequence $\{h(\mu_1^n)\bmod K_v^{\circ\circ}\}_{n\geq0}$ is weakly periodic by Lemma~\ref{algebraic implies weakly periodic}. In particular, looking at the arithmetic progression $(q-1)\mathbb N+b$, there exist integers $k>0$ and $r>r'\geq0$ such that \begin{equation}\label{periodic condition}
    h(\mu_1^{(q-1)(kn+r)+b})= h(\mu_1^{(q-1)(kn+r')+b})\bmod K_v^{\circ\circ}\text{ for all }n\geq0.
\end{equation} 

We denote $r_n=(q-1)(kn+r)+b$ and $r_n'=(q-1)(kn+r')+b$ for convenience. Observe that $r_n-r_n'=(q-1)(r-r')$ for all $n\geq0$. Recall that $\theta=\arg\mu_1\not\in2\pi\mathbb Q$. Thus, the sets $\{\{r_n\theta/2\pi\}:n\in\mathbb N\}$ and $\{\{r'_n\theta/2\pi\}:n\in\mathbb N\}$ are both dense in $(0,1)$. We may assume that $I_1=[\alpha,\beta]$ for some $0\leq \alpha<\beta\leq 1$ with $\beta-\alpha<1$; the case $I_1=[0,1]\setminus(\alpha,\beta)$ is proved similarly. We have two cases.
\begin{itemize}
\item If $\{(q-1)(r-r')\theta/2\pi\}>\beta-\alpha$, then we can find $n\geq0$ such that $\{r_n\theta/2\pi\}\in I_1$ (which is close to the right of $\alpha$) and $\{r'_n\theta/2\pi\}\not\in I_1$. In other words, $\{r_n\theta/2\pi\}\in I_1$ and $\{r_n'\theta/2\pi\}\in I_j$ for some $j\neq 1$. It follows that \[ h(\mu_1^{r_n})=c_1\mu_1^{r_n}+\overline{c_1}\mu_2^{r_n}= c_1\mu_1^b+\overline{c_1}\mu_2^b\bmod K_v^{\circ\circ},\]
and \[ h(\mu_1^{r_n'})=c_j\mu_1^{r_n'}+\overline{c_j}\mu_2^{r_n'}= c_j\mu_1^b+\overline{c_j}\mu_2^b\bmod K_v^{\circ\circ}.\]
Here, we use the fact that $x^{q-1}= 1\bmod K_v^{\circ\circ}$ for every $x\in K_v^{\circ}\setminus K_v^{\circ\circ}$. Since $c_1\mu_1^b+\overline{c_1}\mu_2^b\not= c_j\mu_1^b+\overline{c_j}\mu_2^b\bmod K_v^{\circ\circ}$ by the assumption on $p$, the condition~\eqref{periodic condition} cannot be fulfilled for such $n$, a contradiction.

    \item If $\{(q-1)(r-r')\theta/2\pi\}\leq\beta-\alpha$, then we can find $n\geq0$ such that $\{r_n\theta/2\pi\}\not\in I_1$ (which is close to the left of $\alpha$) and $\{r'_n\theta/2\pi\}\in I_1$. Therefore, $\{r_n\theta/2\pi\}\in I_j$ for some $j\neq 1$ and $\{r_n'\theta/2\pi\}\in I_1$. An argument similar to the previous case also leads to a contradiction.
\end{itemize} 
The proof of Theorem~\ref{dichotomy of the generating power series for toric surfaces} is complete.

\begin{remark}
We cannot avoid the possibility of finitely many exceptional prime numbers $p$ in Theorem~\ref{dichotomy of the generating power series for toric surfaces}. Indeed, $\Delta_{1,\varphi,pD}\bmod p$ is always zero for all prime numbers $p$. A less trivial obstruction comes from the following example. Consider the series $\Delta_{1,\varphi_\mu,\mathcal{O}(1)}(z)$ as in Corollary~\ref{the example of BDJ}. In this case, the degree sequence is given by the values at $A_\mu^n$ of a piecewise linear function $h\colon \Pi\rightarrow\mathbb R$ with coefficients in $\mathbb Z$, see~\cite[Formula~(1.1)]{BDJ-2020}. Therefore, $\Delta_{1,\varphi_\mu,\mathcal{O}(1)}(z)\bmod p$ is equal to $1$ for every common prime divisor $p$ of $\re \mu$ and $\im \mu$.
\end{remark}

\subsection{Proof of Theorem~\ref{when two degree sequences coincide}}

Recall that $\varphi$ and $\varphi'$ are two monomial surface maps which are absolutely irreducible over $\mathbb Q$ with the same degree sequence. We need to show that they are semi-conjugate.

     Let $A$ and $A'$ be their corresponding matrices in $M_2(\mathbb Z)$. Let $\mu_{1},\mu_{2}$ be the eigenvalues of $A$ and $\mu_{1}',\mu_{2}'$ be the eigenvalues of $A'$ ordered so that $|\mu_1|\geq|\mu_2|$ and $|\mu_1'|\geq|\mu_2'|$. Then  $\lambda_1(\varphi)=|\mu_{1}|$ and $\lambda_1(\varphi')=|\mu_{1}'|$. The hypothesis yields $|\mu_1|=\lambda_1(\varphi)=\lambda_1(\varphi')=|\mu_1'|$ and we denote these numbers by $\lambda_1$ for short. We need to show that there exists some $u \in \{1,2,3,4,6\}$ such that $A^u$ and ${A'}^u$ are semi-conjugate via some $P \in M_2(\mathbb{Z})$ with $\det P\neq0$.
     
     Note that since $A$ is absolutely $\mathbb Q$-irreducible, no eigenvalue of $A$ is of the form $x^{1/n}$ for some $x\in\mathbb Q$ and $n\in\mathbb N$. Therefore, either $|\mu_1|>|\mu_2|$ and both eigenvalues are real; or they are complex conjugates with $\mu_1/\mu_2$ not being a root of unity.

     Suppose first that $\mu_1$ and $\mu_2$ are real; then so are $\mu_1'$ and $\mu_2'$. In particular, $A$ and $A'$ are diagonalizable over $\mathbb R$. We note that $\mu_1$ and $\mu_2$ (respectively, $\mu_1'$ and $\mu_2'$) are Galois conjugates.  Thus, it follows from $|\mu_1|=|\mu_1'|$ that $\mu_1^2=\mu_1'^2$ and $\mu_2^2=\mu_2'^2$. In other words, $A^2$ and ${A'}^2$ are matrices with coefficients in $\mathbb Z$ with the same eigenvalues. It then follows from~\cite[A.VII.33 Corollary~3]{Bourbaki-algebraII} that they are semi-conjugate by some $P\in M_2(\mathbb Z)$ with $\det P\neq0$.

       Next, if $\mu_1$ and $\mu_2$ are complex conjugates and $\mu_1/\mu_2$ is not a root of unity, then $\Delta_{1,\varphi,D}(z)$ is not rational, hence neither is $\Delta_{1,\varphi',D'}(z)$. Thus, $\mu_1'$ and $\mu_2'$ are also complex conjugates, and $\mu_1'/\mu_2'$ is not a root of unity. We set $\theta=\arg\mu_1\not\in2\pi\mathbb Q$ and $\theta'=\arg\mu_1'\not\in2\pi\mathbb Q$. Recall that for $z\in\mathbb D(0,\lambda_1^{-1})$, we have 
            $$\Delta_{1,\varphi,D}(z)=\sum_{m\in\mathbb Z}\frac{a_{m}}{1-e^{\irm m\theta}\lambda_1 z}\text{ and }\Delta_{1,\varphi',D'}(z)=\sum_{m\in\mathbb Z}\frac{a_{m}'}{1-e^{\irm m\theta'}\lambda_1 z}.$$
Here, for $m\neq\pm1$, $a_m=\frac{d_m}{m^2-1}$ and $a_m'=\frac{d_m'}{m^2-1}$ for some linear recurrence sequences $\{d_m:m\in\mathbb Z\}$ and $\{d_m':m\in\mathbb Z\}$. Furthermore, both sequences $\{d_m:m\in\mathbb Z\}$ and $\{d_m':m\in\mathbb Z\}$ are non-zero on some arithmetic progressions.
        
Consider the lattice $\Lambda_{\theta,\theta'}\coloneqq \{(m,m'):m\theta= m'\theta'\text{ modulo }2\pi\mathbb Z\}\subset\mathbb Z^2$. Its rank is at most 1 since $\theta$ and $\theta'$ are incommensurable with $2\pi$. Thus it is generated by some $(u_0,u_0')\in\mathbb Z^2$, i.e., $e^{\irm mu_0\theta}=e^{\irm mu_0'\theta'}$ for all $m\in\mathbb Z$. Lemma~\ref{A general result on natural boundary} yields that for $m\in\mathbb Z$, we have
\begin{equation*}\label{eq: singularities of Delta varphi}
  \lim_{\rho\to 1^{-}}(1-\rho)\Delta_{1,\varphi,D}(\rho e^{-\irm mu_0'\theta'})=a_{mu_0}, \text{ and } \lim_{\rho\to 1^{-}}(1-\rho)\Delta_{1,\varphi',D'}(\rho e^{-\irm mu_0\theta}) =a'_{mu_0'}.
\end{equation*}
On the other hand, for every $m_0\in\mathbb Z$ such that $e^{\irm m_0\theta}\not\in\{e^{\irm m'\theta'}:m'\in\mathbb Z\}$, we have
\begin{equation*}\label{eq: vanishing points of (1-rho)Delta varphi'}
\lim_{\rho\to1^{-}}(1-\rho)\Delta_{1,\varphi',D'}(\rho e^{-\irm m_0\theta}\lambda_1)=0.   
\end{equation*}
A similar property holds for $m_0'\in\mathbb Z$ such that $e^{\irm m_0'\theta'}\not\in\{e^{\irm m\theta}:m\in\mathbb Z\}$. 

It then follows from  $\Delta_{1,\varphi,D}(z)=\Delta_{1,\varphi',D'}(z)$ that $\Lambda_{\theta,\theta'}$ is a non-trivial lattice. Therefore, $(u_0,u_0')\neq(0,0)$ and $a_{mu_0}=a'_{mu_0'}$ for all $m\in\mathbb Z$, and they are non-zero for infinitely many $m$.
 It follows that 
 \[\frac{d_{mu_0}}{(mu_0)^2-1}=\frac{d'_{mu_0'}}{(mu_0')^2-1}\] for all $m\neq\pm1$; hence,
 \begin{equation}\label{eq: two linear recurrence sequences coincide}
     m^2(u_0'^2d_{mu_0}-u_0^2d'_{mu_0'})=d_{mu_0}-d'_{mu_0'}\text{ }.
 \end{equation}
 
We will show that $u_0^2=u_0'^2$. Assume the contrary, then $d_{mu_0}\neq d'_{mu_0'}$ for infinitely many $m$, hence $u_0'^2d_{mu_0}\neq u_0^2d'_{mu_0'}$ for infinitely many $m$. It follows that the linear recurrence sequence $\{d_{mu_0}-d'_{mu_0'}:m\in\mathbb Z\}$ is represented by a power sum of degree at most $1$ by Remark~\ref{a closed form for dm} while the linear recurrence sequence
$\{m^2(u_0'^2d_{mu_0}-u_0^2d'_{mu_0'}):m\in\mathbb Z\}$ is represented by a power sum of degree at least $2$, so the equality~\eqref{eq: two linear recurrence sequences coincide} is absurd.

        Therefore $u_0^2=u_0'^2$. We then set $u=|u_0|$, whence either $u(\theta-\theta')\in 2\pi\mathbb Z$ or $u(\theta+\theta')\in 2\pi\mathbb Z$. This means that $e^{\irm u\theta}=e^{\pm\irm u\theta'}$, hence
        $A^{u}$ and $A'^{u}$ are matrices with coefficients in $\mathbb Z$ and with the same eigenvalues. Therefore, they are semi-conjugate by some $P\in M_2(\mathbb Z)$ with $\det P\neq0$ by~\cite[A.VII.33 Corollary~3]{Bourbaki-algebraII}. 

        It remains to show $u\in\{1,2,3,4,6\}$. Observe that $\xi\coloneqq e^{\irm (\theta-\theta')}$ is a primitive $u$-th root of unity, hence $\mu_1=\xi\mu_1'$. Since $\mu_1$ and $\mu_1'$ are quadratic integers, $\mathbb Q(\mu_1,\mu_1')$ is a Galois extension of $\mathbb Q$ of degree $[\mathbb Q(\mu_1,\mu_1'):\mathbb Q]\leq4$. Since $\mathbb Q(\xi)\subseteq\mathbb Q(\mu_1,\mu_1')$ and $[\mathbb Q(\xi):\mathbb Q]=\phi(u)$, we have $\phi(u)|4$. Here, $\phi$ is Euler's totient function. Thus $u\in\{1,2,3,4,5,6,8,10,12\}$. 
        
        Next, assume that $u\in\{5,8,10,12\}$, so that $\phi(u)=4$. Thus $[\mathbb Q(\mu_1,\mu_1'):\mathbb Q]=[\mathbb Q(\xi):\mathbb Q]=4$, whence $\mathrm{Gal}(\mathbb Q(\mu_1,\mu_1')/\mathbb Q)\cong\mathbb Z/2\mathbb Z\times\mathbb Z/2\mathbb Z$. Therefore $u\neq 5,10$ since both  $\mathrm{Gal}(\mathbb Q(e^{\irm 2\pi/5})/\mathbb Q)$ and $\mathrm{Gal}(\mathbb Q(e^{\irm 2\pi/10})/\mathbb Q)$ are cyclic of order $4$. Thus, $u\in\{8,12\}$. We consider the case $u=12$; the case $u=8$ is proved similarly. In this case, $\xi$ is a primitive 12th root of unity. It is known that the minimal polynomial of $\xi$ is $X^4-X^2+1$, and $\mathbb Q(\xi)$ contains only three quadratic subfields: $\mathbb Q(\sqrt{3})$, $\mathbb Q(\irm)$, and $\mathbb Q(\sqrt{-3})$. Since $\mathbb Q(\mu_1)$ and $\mathbb Q(\mu_1')$ are two distinct complex quadratic subfields of $\mathbb Q(\xi)$, one of them must be $\mathbb Q(\irm)$, and the other is $\mathbb Q(\sqrt{-3})$. We may assume that $\mathbb Q(\mu_1)=\mathbb Q(\irm)=\mathbb Q(\xi^3)$ and $\mathbb Q(\mu_1')=\mathbb Q(\sqrt{-3})=\mathbb Q(\xi^4)$. Thus there are $a,b,c,d\in\mathbb Q$ such that $\mu_1=a+b\xi^3$ and $\mu_1'=c+d\xi^4=c-d+d\xi^2$. From $\mu_1=\xi\mu_1'$, we obtain
        $$a+b\xi^3=\xi(c-d+d\xi^2)=(c-d)\xi+d\xi^3.$$
        Thus $a=0$ and $b=c=d$; hence, $\mu_1=\pm b\irm$, which contradicts the hypothesis $\arg\mu_1\not\in2\pi\mathbb Q$.
        
        We conclude that $u\not\in\{5,8,10,12\}$. The theorem is proved.

\begin{remark}
    Observe that the irreducibility assumption is necessary. Indeed, consider monomial self-maps $\varphi$ and $\varphi'$ of $\mathbb P^2$ associated with matrices \[A=\begin{pmatrix}
        a&0\\0&b
    \end{pmatrix}\text{ and }A'=\begin{pmatrix}
        a&0\\0&c
    \end{pmatrix}\]
    where $a,b,c$ are positive integers and $a>b>c$. Let $D=D'=\mathcal{O}(1)$. Then \[\Delta_{1,\varphi,\mathcal{O}(1)}(z)=\Delta_{1,\varphi',\mathcal{O}(1)}(z)=\sum_{n\geq0}a^nz^n=\frac{1}{1-az}.\] However, $A^{u}$ and $A'^{u}$ are not conjugate for any $u>0$.
\end{remark}

\bibliographystyle{plain}

\begin{thebibliography}{99}

\bibitem{Abarenkova-Auriac-Boukraa-Maillard-1999} N. Abarenkova, J.-Ch. {Anglès d’Auriac}, S. Boukraa, and J.-M. Maillard. \textit{Growth-complexity spectrum of some discrete dynamical systems}. Physica D \textbf{130}(1), 27--42 (1999)

\bibitem{Adamczewski-Bell-2013} B. Adamczewski, and J. P. Bell. \textit{Diagonalization and rationalization of algebraic Laurent series}. Ann. Sci. \'Ec. Norm. Sup\'er. \textbf{46}(6), 963--1004 (2013)

\bibitem{Adamczewski-Bell-Delaygue-2019} B. Adamczewski, J. P. Bell, and \'E. Delaygue. \textit{Algebraic independence of $G$-functions and congruences ``\`a la Lucas''}. Ann. Sci. \'Ec. Norm. Sup\'er. \textbf{52}(3), 515--559 (2019)

\bibitem{Adamczewski-Bostan-Caruso-2023} B. Adamczewski, A. Bostan, and X. Caruso. \textit{Diagonals and algebraicity modulo $p$: a sharper degree bound}. To appear in Ann. Sci. \'Ec. Norm. Sup\'er. \texttt{arXiv:2601.14920}


\bibitem{Adamczewski-Dreyfus-Hardouin-2021} B. Adamczewski, T. Dreyfus, and C. Hardouin. \textit{Hypertranscendence and linear difference equations}. J. Amer. Math. Soc. \textbf{34}(2), 475--503 (2021)

\bibitem{Becker-Topfer-94} P. Becker, and T. T\"{o}pfer. \textit{Transcendency Results for Sums of Reciprocals of Linear Recurrences}. Math. Nachr. \textbf{168}, 5--17 (1994)

\bibitem{Bedford-Kim-2004} E. Bedford, and K. Kim. \textit{On the degree growth of birational mappings in higher dimension}. J. Geom. Anal. \textbf{14}(4), 567--596 (2004)

\bibitem{Bell-talk-NTWS} J. P. Bell. \textit{A transcendental dynamical degree (Number Theory Web Seminar 052)}. 2024.

\bibitem{Bell-Coons-Rowland-2013} J. P. Bell, M. Coons, and E. Rowland. \textit{The rational-transcendental dichotomy of Mahler functions}. J. Integer Seq. \textbf{16}(2), article 13.2.10 (2013)


\bibitem{BDJ-2020} J. P. Bell, J. Diller, and M. Jonsson. \textit{A transcendental dynamical degree}. Acta Math. \textbf{225}(2), 193--225 (2020)

\bibitem{BDJK-2024} J. P. Bell, J. Diller, M. Jonsson, and H. Krieger. \textit{Birational maps with transcendental dynamical degree}. Proc. Lond. Math. Soc. \textbf{128}(1), e12573 (2024)



\bibitem{Nguyen-Bell-Gunn-Saunders-2023} J. P. Bell, K. Gunn, K. Nguyen, and J. C. Saunders. \textit{A general criterion for the P\'olya-Carlson dichotomy and application}. Trans. Amer. Math. Soc. \textbf{376}(6), 4361--4382 (2023)

\bibitem{Boucksom-Favre-Jonsson-2008} S. Boucksom, C. Favre, and M. Jonsson. \textit{Degree growth of meromorphic surface maps}. Duke Math. J. \textbf{141}(3), 519--538 (2008)




\bibitem{Boukraa-Hassani-Maillard-2003} S. Boukraa, S. Hassani, and J.-M. Maillard. \textit{Noetherian mappings}. Physica D \textbf{185}(1), 3--44 (2003)

\bibitem{Boudreau-Holmes-Nguyen-2025} F. B. Boudreau, E. Holmes, and D.-K. Nguyen. \textit{Adelic perturbation of rational functions and applications}. Math. Ann. \textbf{392}(2), 2253--2275 (2025) 



\bibitem{Bourbaki-algebraII} N. Bourbaki. \textit{Algebra II: Chapters 4--7.} Elements of Mathematics (Berlin) Translated from the French by P. M. Cohn and J. Howie. Springer-Verlag, Berlin (1990).



\bibitem{Bousquet-Melou-2006} M. Bousquet-M\'elou. \textit{Rational and algebraic series in combinatorial enumeration}. In International Congress of Mathematicians III, Eur. Math. Soc., 2006, 789--826.

\bibitem{Byszewski-Konieczny-2020} J. Byszewski, and J. Konieczny. \textit{Automatic Sequences and Generalised Polynomials}. Canad. J. Math. \textbf{72}(2), 392--426 (2020)





\bibitem{Caruso-Furnsinn-Vargas-Montoya-2025} X. Caruso, F. F\"{u}rnsinn, and D. Vargas-Montoya. \textit{Galois groups of reductions modulo $p$ of $D$-finite series}. Preprint (2025) \texttt{arXiv:2504.09429}

\bibitem{Christol-1979} G. Christol. \textit{Ensembles presque periodiques $k$-reconnaissables}. Theor. Comput. Sci. \textbf{9}(1), 141--145 (1979)

\bibitem{Christol-Kamae-France-Rauzy-1980} G. Christol, T. Kamae, M. Mend\`es France, and G. Rauzy.  \textit{Suites alg\'ebriques, automates et substitutions}. Bull. Soc. Math. France \textbf{108}, 401--419 (1980)

\bibitem{Dang-2020} N.-B. Dang. \textit{Degrees of iterates of rational maps on normal projective varieties}. Proc. Lond. Math. Soc. \textbf{121}(5), 1268--1310 (2020)

\bibitem{Dang-Ramadas-2021} N.-B. Dang, and R. Ramadas. \textit{Dynamical invariants of monomial correspondences}. Ergodic Theory Dynam. Syst. \textbf{41}(7), 2000--2015 (2021)


\bibitem{Deligne-1984} P. Deligne. \textit{Int\'egration sur un cycle \'evanescent}. Invent. Math. \textbf{76}(1), 129--143 (1984)

\bibitem{Diller-Favre-2001} J. Diller, and C. Favre. \textit{Dynamics of bimeromorphic maps of surfaces}. Amer. J. Math. \textbf{123}(6), 1135--1169 (2001)

\bibitem{Dinh-Sibony-2005} T.-C. Dinh et N. Sibony. \textit{Une borne sup\'erieure pour l’entropie topologique d’une application rationnelle}. Ann. of Math. \textbf{161}(3), 1637--1644 (2005)

\bibitem{Favre-2003} C. Favre. \textit{Les applications monomiales en deux dimensions}. Michigan Math. J. \textbf{51}(3), 467--475 (2003)


\bibitem{Favre-Wulcan-2012} C. Favre, and E. Wulcan. \textit{Degree growth of monomial maps and McMullen’s polytope algebra}. Indiana Univ. Math. J. \textbf{61}(2), 493--524 (2012)

\bibitem{Fornaess-Sibony-1995} J. E. Forn{\ae}ss, and N. Sibony. \textit{Complex dynamics in higher dimension. II}. In: Bloom, T., Catlin, D.W., D'Angelo, J.P., Siu, Y.-T. (eds.) Modern Methods in Complex Analysis: The Princeton Conference in Honor of Gunning and Kohn. Annals of Mathematics Studies, vol. 137, pp. 135--182. Princeton University Press, Princeton (1995)

\bibitem{Hasselblatt-Propp-2007} B. Hasselblatt, and J. Propp. \textit{Degree-growth of monomial maps}. Ergodic Theory Dynam. Syst. \textbf{27}(5), 1375--1397 (2007)


\bibitem{Hilbert-1902} D. Hilbert. \textit{Mathematical problems}. Bull. Amer. Math. Soc. \textbf{8}(10), 437--479 (1902)

\bibitem{Hormander-book-1983} L. H\"{o}rmander. \textit{The Analysis of Linear Partial Differential Operators I: Distribution Theory and Fourier Analysis.} Grundlehren der mathematischen Wissenschaften, vol. 256. Springer, Berlin, Heidelberg (1983).


\bibitem{Huber-Sturmfels-1995} B. Huber, and B. Sturmfels. \textit{A Polyhedral Method for Solving Sparse Polynomial Systems}. Math. Comput. \textbf{64}(212), 1541--1555 (1995)


\bibitem{Ji-Xie-2023} Z. Ji, and J. Xie. \textit{Local rigidity of Julia sets}.  to appear in Amer. J. Math. \texttt{arXiv:2302.02562}


\bibitem{Jonsson-Wulcan-2011} M. Jonsson, and E. Wulcan. \textit{Stabilization of monomial maps}. Michigan Math. J. \textbf{60}(3), 629--660 (2011)



\bibitem{Kontsevich-Zagier-2001} M. Kontsevich, and D. Zagier. \textit{Periods}. In: Engquist, B., Schmid, W. (eds.) Mathematics Unlimited---2001 and Beyond, pp. 771--808. Springer, Berlin, Heidelberg (2001)

\bibitem{Lin-2012-BSMF} J.-L. Lin. \textit{Pulling back cohomology classes and dynamical degrees of monomial maps}. Bull. Soc. Math. France \textbf{140}(4), 533--549 (2012)

\bibitem{Lin-Wulcan-2014} J.-L. Lin, and E. Wulcan. \textit{Stabilization of monomial maps in higher codimension}. Ann. Inst. Fourier \textbf{64}(5), 2127--2146 (2014)

\bibitem{Milnor-Thurston-1988} J. Milnor, and W. Thurston. \textit{On iterated maps of the interval}. In: Alexander, J.C. (ed.) Dynamical Systems. Lecture Notes in Mathematics, vol. 1342, pp. 465--563. Springer, Berlin, Heidelberg (1988)

\bibitem{Oda-toric-geometry} T. Oda. \textit{Convex bodies and algebraic geometry: An introduction to the theory of toric varieties}. Ergebnisse der Mathematik und ihrer Grenzgebiete (3), vol. 15. Translated from the Japanese. Springer-Verlag, Berlin, viii+212 pp. (1988)

\bibitem{Schneider-book-2013} R. Schneider. \textit{Convex Bodies: The Brunn--Minkowski Theory, 2nd edn.} Encyclopedia of Mathematics and its Applications. Cambridge University Press (2013).


\bibitem{Stanley-1980} R. P. Stanley. \textit{Differentiably Finite Power Series}. Eur. J. Combin. \textbf{1}(2), 175--188 (1980)


\bibitem{Stanley-Enumerative-Combinatorics-2011} R. P. Stanley. \textit{Enumerative Combinatorics, vol. 1, 2nd edn.} Cambridge Studies in Advanced Mathematics. Cambridge University Press (2012)

\bibitem{Stanley-Fomin-vol2-1999} R. P. Stanley, and S. Fomin. \textit{Enumerative Combinatorics: vol. 2}. Cambridge Studies in Advanced Mathematics, vol. 62. Cambridge University Press, Cambridge, 1999.


\bibitem{Sugimoto-2025} Y. Sugimoto. \textit{A birational map of a projective space whose intermediate dynamical degrees are all transcendental}. Preprint (2025) \texttt{arXiv:2503.00688}

\bibitem{Truong-2020} T.-T. Truong. \textit{Relative dynamical degrees of correspondences over a field of arbitrary characteristic}. J. Reine Angew. Math. \textbf{758}, 139--182 (2020)

\bibitem{Urech-2018} C. Urech. \textit{Remarks on the degree growth of birational transformations}. Math. Res. Lett. \textbf{25}(1), 291--308 (2018)

\bibitem{Vargas-Montoya-2021} D. Vargas-Montoya. \textit{Alg\'ebricit\'e modulo $p$, s\'eries hyperg\'eom\'etriques et structures de Frobenius fortes}. Bull. Soc. Math. France \textbf{149}(3), 439--477 (2021)
\end{thebibliography}

{Institut Camille Jordan, Universit\'e Claude Bernard Lyon 1\\ 21 avenue Claude Bernard, 69100 Villeurbanne, France}\\Email: \texttt{nguyen@math.univ-lyon1.fr} 
\end{document}